\pdfoutput=1
\documentclass[11pt]{article}
\usepackage[left=1in,right=1in,top=1in,bottom=1in]{geometry}
\usepackage{times}
\usepackage{expl3}
\usepackage{cite}
\usepackage[table]{xcolor}
\usepackage{multirow}
\usepackage{stackengine} 
\usepackage{hhline}
\usepackage{lipsum}
\usepackage{titlesec}
\usepackage{wrapfig}
\usepackage{enumerate}
\usepackage{epsfig}
\usepackage{amsmath}
\usepackage{tabularx}
\usepackage{array}
\usepackage{booktabs}
\usepackage{enumitem}
\usepackage{bbm}
\usepackage{calc}
\usepackage{graphicx}
\usepackage{amsmath}
\usepackage[title]{appendix}
\usepackage{amssymb}
\usepackage{epstopdf}
\usepackage{boldline}
\usepackage{arydshln}
\usepackage{calligra}
\usepackage{bm}
\usepackage{url}
\usepackage{blindtext}
\usepackage{accents}

\newcommand{\define}{\stackrel{\mbox{\tiny def}}{=}}

\newtheorem{definition}{Definition}
\newtheorem{theorem}{Theorem}
\newtheorem{corollary}{Corollary}
\newtheorem{lemma}{Lemma}

\usepackage{mathtools}
\usepackage{epstopdf}
\usepackage{balance}
\usepackage{thmtools}
\usepackage{thm-restate}
\usepackage{hyperref}
\usepackage{cleveref}
\usepackage[mathscr]{euscript}

\usepackage[ruled,vlined]{algorithm2e}
\include{pythonlisting}
\newcommand{\ostar}{\mathbin{\mathpalette\make@circled\star}}

\makeatletter
\newcommand{\removelatexerror}{\let\@latex@error\@gobble}
\makeatother
\makeatletter
\newcommand*{\rom}[1]{\expandafter\@slowromancap\romannumeral #1@}
\makeatother

\ExplSyntaxOn
\newcommand\latinabbrev[1]{
  \peek_meaning:NTF . {
    #1\@}%
  { \peek_catcode:NTF a {
      #1.\@ }%
    {#1.\@}}}
\ExplSyntaxOff

\titleclass{\subsubsubsection}{straight}[\subsubsection]

\begin{document}
\vspace{1cm}
\title{Multivariate Mond-Pecaric Method with Applications to Hypercomplex Function Sobolev Embedding}\vspace{1.8cm}
\author{Shih-Yu~Chang
\thanks{Shih Yu Chang is with the Department of Applied Data Science,
San Jose State University, San Jose, CA, U. S. A. (e-mail: {\tt
shihyu.chang@sjsu.edu}). 
           }}

\maketitle

\begin{abstract}
Mond and Pecaric introduced a method to simplify the determination of complementary inequalities for Jensen's inequality by converting it into a single-variable maximization or minimization problem of continuous functions. This principle has significantly enriched the field of operator inequalities. Our contribution lies in extending the Mond-Pecaric method from single-input operators to multiple-input operators. We commence by defining normalized positive linear maps, accompanied by illustrative examples. Subsequently, we employ the Mond-Pecaric method to derive fundamental inequalities for multivariate hypercomplex functions bounded by linear functions. These foundational inequalities serve as the basis for establishing several multivariate hypercomplex function inequalities, focusing on ratio relationships. Additionally, we present similar results based on difference relationships. Finally, we apply the derived multivariate hypercomplex function inequalities to establish Sobolev embedding via Sobolev inequality for hypercomplex functions with operator inputs.
\end{abstract}


\section{Introduction}\label{sec: Introduction}

The Kantorovich inequality states that when $\bm{A}$ is a positive operator on the Hilbert space $\mathfrak{H}$ and satisfies the conditions $m \bm{I}_{\mathfrak{H}} \leq \bm{A} \leq M \bm{I}_{\mathfrak{H}}$ for positive scalars $m, M$ with $m < M$ and the identity operator $\bm{I}_{\mathfrak{H}}$ in the Hilber space $\mathfrak{H}$, then
\begin{eqnarray}
\langle\bm{A}\bm{x},\bm{x}\rangle \langle\bm{A}^{-1}\bm{x},\bm{x}\rangle \leq \frac{(M+m)^2}{4Mm},
\end{eqnarray}
where $\bm{x} \in \mathfrak{H}$ with $\left\|\bm{x}\right\|=1$. Several authors have explored extensions of the Kantorovich inequality. Notably, Mond and Pecaric have conducted an extensive research series in this area~\cite{furuta2005mond,fujii2012recentMP}. They developed a method that generates complementary inequalities to Jensen's type inequalities associated with convex functions. Consequently, they provided complementary inequalities to the Holder-McCarthy inequality and extensions of the Kantorovich type. Furuta~\cite{furuta1995extension} further extended the work by Ky Fan~\cite{fan1966some} and the Mond-Pecaric generalizations of the Kantorovich inequality by incorporating the ideas of Ky Fan, Mond, and Pecaric. Additionally, in integral expressions, S.-E. Takahasi et al.~\cite{takahashi1999inverse} introduced another formula for a complementary inequality to Jensen's inequality, which includes the Kantorovich inequality as a special case. By amalgamating the concepts of Furuta and Takahasi, new insights into the method established by Mond and Pecaric are revealed, allowing for the derivation of complementary inequalities to Jensen's inequality for convex functions.

Mond and Pecaric proposed a method to transform the complementary inequality determination of Jensen's inequality to a single-variable maximization or minimization problem of continuous functions. Employing this approach, we can address general complementary inequalities to Jensen's inequalities for convex functions. This formulation simplifies the concept of complementary inequalities and enhances clarity in notions and proofs. This viewpoint proves valuable for exploring topics such as the Hadamard product, positive linear maps, operator means, and order-preserving operator inequalities. Referred to as the Mond-Pecaric method, this principle yields substantial results in the realm of operator inequalities~\cite{furuta2005mond,fujii2012recentMP}. For instance, the operator Jensen's inequality extends the concept of convexity to operator variables, offering insights into function convexity concerning matrices or tensors, which are commonly encountered in various mathematical and engineering tasks \cite{chang2024TZF, chang2024generalizedCDJ, chang2023TenLevenberg, chang2023BMSBSVD, chang2023personalized, chang2023TenLeastSquares, chang2023TenProj, chang2022TLMS, chang2021TenInv, chang2022TWF, chang2022Kalman, chang2022TCASI, chang2021RLS}.

The contribution of this work is to extend the Mond-Pecaric method from single input operator to multiple input operators. We begin by presenting the definition of normalized positive linear maps along with illustrative examples. Then, the Mond-Pečarić method is utilized to derive fundamental inequalities for multivariate hypercomplex functions bounded by linear functions. These fundamental inequalities are used to establish several multivariate hypercomplex function inequalities based on ratio relationships. Besides, we also present analogous results but based on difference relationships. Finally, we apply the derived multivariate hypercomplex function inequalities to establish Sobolev inequality for hypercomplex functions with operator inputs.

The definition about the normalized positive linear map and its examples are provided by Section~\ref{sec: Normalized Positive linear maps}. Fundamental inequalities for multivariate hypercomplex functions bounded by linear functions are derived by Mond-Pecaric method in Section~\ref{sec: Fundamental Inequalities for Multivariate Hypercomplex Functions}. In Section~\ref{sec: Multivariate Hypercomplex Functions Inequalities: Ratio Kind}, several multivariate hypercomplex functions inequalities are established based on ratio relationship. On the other hand, in Section~\ref{sec: Multivariate Hypercomplex Functions Inequalities: Difference Kind}, several multivariate hypercomplex functions inequalities are established based on difference relationships. Finally, the application of derived multivariate hypercomplex function inequalities to established Sobolev inequality for hypercomplex function with operator inputs is discussed in  Section~\ref{sec: Applications}. 

\textbf{Nomenclature:} 
Inequalities $\geq, >, \leq,$ and $<$, when applied to operators, follow the Loewner order. The symbol $\Lambda(\bm{A})$ denotes the spectrum of the operator $\bm{A}$, i.e., the set of eigenvalues of $\bm{A}$. If $\Lambda(\bm{A})$ consists of real numbers, $\min(\Lambda(\bm{A}))$ and $\max(\Lambda(\bm{A}))$ represent the minimum and maximum values within $\Lambda(\bm{A})$, respectively. In this work, a hypercomplex function indicates a function with operators/matrices/tensors as input variables.  We assume that all input operators of multivariate functions discussed in this work are self-adjoint, and the output operator is self-adjoint too. The summation of polynomials with respect to each self-adjoint operator is a class of multivariate functions with self-adjoint operators as their outputs with self-adjoint input operators. 

\section{Normalized Positive Linear Maps}\label{sec: Normalized Positive linear maps}

In this section, we will introduce the notion about \emph{Normalized Positive Linear Maps}. Let us set up several symbols for later usage. $\mathscr{B}(\mathfrak{H})$ represents semi-algebra of all bounded linear operators on a Hilbert space from $\mathfrak{H}$ to $\mathfrak{H}$. $\bm{I}_{\mathfrak{H}}$ is the identity operator in  $\mathfrak{H}$.  $\mathscr{B}_{sa}(\mathfrak{H})$ represents semi-algebra of all self-adjoint bounded linear operators on a Hilbert space from $\mathfrak{H}$ to $\mathfrak{H}$.

In the following definition~\ref{def: normalized positive linear map}, we define a normalized positive linear map~\cite{furuta2005mond,fujii2012recentMP}.
\begin{definition}\label{def: normalized positive linear map}
A map $\Phi: \mathscr{B}(\mathfrak{H}) \rightarrow \mathscr{B}(\mathfrak{K})$ is considered a normalized positive linear map if it satisfies the following conditions:
\begin{itemize}
\item Linearity: $\Phi(a\bm{X}+b\bm{Y})=a\Phi(\bm{X})+b\Phi(\bm{Y})$ for any $a,b \in \mathbb{C}$ and any $\bm{X}, \bm{Y} \in \mathscr{B}(\mathfrak{H})$.
\item Positivity: If $\bm{X}\geq\bm{Y}$, then $\Phi(\bm{X})\geq\Phi(\bm{Y})$.
\item Normalization: $\Phi(\bm{I}_{\mathfrak{H}})=\bm{I}_{\mathfrak{K}}$, where $\bm{I}_{\mathfrak{H}}$ and $\bm{I}_{\mathfrak{K}}$ are the identity operators of the Hilbert spaces $\mathfrak{H}$ and $\mathfrak{K}$, respectively. 
\end{itemize}
\end{definition}

Given a set of operators $\bm{V}_i \in \mathscr{B}(\mathfrak{H})$ for $i=1,2,\ldots,n$ satisfying 
\begin{eqnarray}
\sum\limits_{i=1}^n \bm{V}_i^{\ast}\bm{V}_i&=&\bm{I}_{\mathfrak{H}},
\end{eqnarray}
where $^{\ast}$ represents the self-adjoint operation. Then the map $\Phi: \mathscr{B}(\mathfrak{H}) \rightarrow \mathscr{B}(\mathfrak{H})$ defined below is a normalized positive linear map.
\begin{eqnarray}
\Phi(\bm{X})&=&\sum\limits_{i=1}^n \bm{V}_i^{\ast}\bm{X}\bm{V}_i. 
\end{eqnarray}

\section{Fundamental Inequalities for Multivariate Hypercomplex Functions}\label{sec: Fundamental Inequalities for Multivariate Hypercomplex Functions}

In this section, we will establish two fundamental inequalities for multivariate hypercomplex functions. 

\begin{theorem}\label{thm: main 2.3}
Let $\bm{A}_{j_i}$ be self-adjoint operators with $\Lambda(\bm{A}_{j_i}) \in [m_i, M_i]$ for real scalars $m_i <  M_i$. The mappings $\Phi_{j_1,\ldots,j_n}: \mathscr{B}(\mathfrak{H}) \rightarrow \mathscr{B}(\mathfrak{K})$ are normalized positive linear maps, where $j_i=1,2,\ldots,k_i$ for $i=1,2,\ldots,n$. We have $n$ probability vectors $\bm{w}_i =[w_{i,1},w_{i,2},\cdots, w_{i,k_i}]$ with the dimension $k_i$ for $i=1,2,\ldots,n$, i.e., $\sum\limits_{\ell=1}^{k_i}w_{i,\ell} = 1$. Let $f(x_1,x_2,\ldots,x_n)$ be any real valued continuous functions with $n$ variables defined on the range $\bigtimes\limits_{i=1}^n [m_i, M_i] \in \mathbb{R}^n$, where $\times$ is the Cartesian product. Besides, we assume that the function $f$ satisfies the following:
\begin{eqnarray}\label{eq f bounds: thm:main 2.3}
\sum\limits_{i=1}^n a_i x_i + b \leq f(x_1,x_2,\ldots,x_n) \leq \sum\limits_{i=1}^n c_i x_i + d,
\end{eqnarray}
for $[x_1,x_2,\ldots,x_n] \in \bigtimes\limits_{i=1}^n [m_i,M_i]$. The function $g(x_1,x_2,\ldots,x_n)$ is also a real valued continuous function with $n$ variables defined on the range $\bigtimes\limits_{i=1}^n [m_i, M_i]$. We also have a real valued function $F(u,v)$ with operator monotone on the first variable $u$ defined on $U \times V$ such that $f(x_1,x_2,\ldots,x_n) \subset U$, and $g(x_1,x_2,\ldots,x_n) \subset V$. Then, we have the following upper bound:
\begin{eqnarray}\label{eq UB: thm:main 2.3}
F\Bigg(\sum\limits_{j_{1}=1,\ldots,j_{n}=1}^{k_1,\ldots,k_n}  w_{j_1}\ldots w_{j_n}\Phi_{j_1,\ldots,j_n}(f(\bm{A}_{j_1},\bm{A}_{j_2},\ldots,\bm{A}_{j_n})),~~~~~~~~~~~~~~~~~~~~~~~~~~~~~~~~~~~~~~~~~~~~~~~~~~\nonumber \\
g\Bigg(\sum\limits_{j_{1}=1,\ldots,j_{n}=1}^{k_1,\ldots,k_n}w_{j_1}\ldots w_{j_n}\Phi_{j_1,\ldots,j_n}(\bm{A}_{j_1}),\ldots,\sum\limits_{j_{1}=1,\ldots,j_{n}=1}^{k_1,\ldots,k_n}w_{j_1}\ldots w_{j_n}\Phi_{j_1,\ldots,j_n}(\bm{A}_{j_n})\Bigg)
\Bigg) \nonumber \\
\leq
\max\limits_{m_i \leq x_i \leq M_i}F\Big(\sum\limits_{i=1}^n c_i x_i + d, g(x_1,x_2,\ldots,x_n)\Big)\bm{I}_{\mathfrak{K}}.
\end{eqnarray}
Similarly, we also have the following lower bound:
\begin{eqnarray}\label{eq LB: thm:main 2.3}
F\Bigg(\sum\limits_{j_{1}=1,\ldots,j_{n}=1}^{k_1,\ldots,k_n}  w_{j_1}\ldots w_{j_n}\Phi_{j_1,\ldots,j_n}(f(\bm{A}_{j_1},\bm{A}_{j_2},\ldots,\bm{A}_{j_n})),~~~~~~~~~~~~~~~~~~~~~~~~~~~~~~~~~~~~~~~~~~~~~~~~~~\nonumber \\
g\Bigg(\sum\limits_{j_{1}=1,\ldots,j_{n}=1}^{k_1,\ldots,k_n}w_{j_1}\ldots w_{j_n}\Phi_{j_1,\ldots,j_n}(\bm{A}_{j_1}),\ldots,\sum\limits_{j_{1}=1,\ldots,j_{n}=1}^{k_1,\ldots,k_n}w_{j_1}\ldots w_{j_n}\Phi_{j_1,\ldots,j_n}(\bm{A}_{j_n})\Bigg)
\Bigg) \nonumber \\
\geq
\min\limits_{m_i \leq x_i \leq M_i}F\Big(\sum\limits_{i=1}^n a_i x_i + b, g(x_1,x_2,\ldots,x_n)\Big)\bm{I}_{\mathfrak{K}}.
\end{eqnarray}
\end{theorem}
\textbf{Proof:}
We will prove Eq.~\eqref{eq UB: thm:main 2.3} first. Due to $f(x_1,x_2,\ldots,x_n) \leq \sum\limits_{i=1}^n c_i x_i + d$ for $[x_1,x_2,\ldots,x_n] \in \bigtimes\limits_{i=1}^n [m_i,M_i]$ from Eq.~\eqref{eq f bounds: thm:main 2.3}, we have
\begin{eqnarray}\label{eq1 UB: thm:main 2.3}
f(\bm{A}_{j_1},\bm{A}_{j_2},\ldots,\bm{A}_{j_n}) \leq \sum\limits_{i=1}^n c_i \bm{A}_{j_i} + d\bm{I}_{\mathfrak{H}},
\end{eqnarray}
where $j_i=1,2,\ldots,k_i$ for $i=1,2,\ldots,n$. Since $\Phi_{j_1,\ldots,j_n}$ is a normalized positive linear map,  we also have
\begin{eqnarray}\label{eq2 UB: thm:main 2.3}
\Phi_{j_1,\ldots,j_n}(f(\bm{A}_{j_1},\bm{A}_{j_2},\ldots,\bm{A}_{j_n}))&\leq&\Phi_{j_1,\ldots,j_n}\Big(\sum\limits_{i=1}^n c_i \bm{A}_{j_i} + d\bm{I}_{\mathfrak{H}}\Big)\nonumber \\
&=& \sum\limits_{i=1}^n c_i \Phi_{j_1,\ldots,j_n}(\bm{A}_{j_i}) +  d\bm{I}_{\mathfrak{K}},
\end{eqnarray}
where $j_i=1,2,\ldots,k_i$ for $i=1,2,\ldots,n$. By multiplying $w_{j_1}\ldots w_{j_n}$ to both sides of Eq.~\eqref{eq2 UB: thm:main 2.3} and summing items with respect to $w_{j_1}\ldots w_{j_n}$, we have 
\begin{eqnarray}\label{eq3 UB: thm:main 2.3}
\lefteqn{\sum\limits_{j_{1}=1,\ldots,j_{n}=1}^{k_1,\ldots,k_n}  w_{j_1}\ldots w_{j_n}\Phi_{j_1,\ldots,j_n}(f(\bm{A}_{j_1},\bm{A}_{j_2},\ldots,\bm{A}_{j_n}))}\nonumber \\
&\leq&
\sum\limits_{i=1}^n c_i \sum\limits_{j_{1}=1,\ldots,j_{n}=1}^{k_1,\ldots,k_n}w_{j_1}\ldots w_{j_n}\Phi_{j_1,\ldots,j_n}(\bm{A}_{j_i})+d\bm{I}_{\mathfrak{K}}.
\end{eqnarray}
Because $m_i \bm{I}_{\mathfrak{H}} \leq \bm{A}_{j_i} \leq M_i \bm{I}_{\mathfrak{H}}$ where $j_i=1,2,\ldots,k_i$ for $i=1,2,\ldots,n$, we have 
\begin{eqnarray}\label{eq4 UB: thm:main 2.3}
\Lambda\left(\sum\limits_{j_{1}=1,\ldots,j_{n}=1}^{k_1,\ldots,k_n}w_{j_1}\ldots w_{j_n}\Phi_{j_1,\ldots,j_n}(\bm{A}_{j_i})\right) \in [m_i, M_i].
\end{eqnarray}
By utilizing the operator monotonicity of $F(u,v)$ with respect to the variable $u$, we have
\begin{eqnarray}\label{eq5 UB: thm:main 2.3}
F\Bigg(\sum\limits_{j_{1}=1,\ldots,j_{n}=1}^{k_1,\ldots,k_n}  w_{j_1}\ldots w_{j_n}\Phi_{j_1,\ldots,j_n}(f(\bm{A}_{j_1},\bm{A}_{j_2},\ldots,\bm{A}_{j_n})),~~~~~~~~~~~~~~~~~~~~~~~~~~~~~~~~~~~~~~~~~~~~~~~~~~\nonumber \\
g\Bigg(\sum\limits_{j_{1}=1,\ldots,j_{n}=1}^{k_1,\ldots,k_n}w_{j_1}\ldots w_{j_n}\Phi_{j_1,\ldots,j_n}(\bm{A}_{j_1}),\ldots,\sum\limits_{j_{1}=1,\ldots,j_{n}=1}^{k_1,\ldots,k_n}w_{j_1}\ldots w_{j_n}\Phi_{j_1,\ldots,j_n}(\bm{A}_{j_n})\Bigg)
\Bigg)\nonumber \\
\leq_1 F\Bigg(\sum\limits_{i=1}^n c_i \sum\limits_{j_{1}=1,\ldots,j_{n}=1}^{k_1,\ldots,k_n}w_{j_1}\ldots w_{j_n}\Phi_{j_1,\ldots,j_n}(\bm{A}_{j_i})+d\bm{I}_{\mathfrak{K}},~~~~~~~~~~~~~~~~~~~~~~~~~~~~~~~~~~~~~~~~~~~~~~~~~~\nonumber \\
g\Bigg(\sum\limits_{j_{1}=1,\ldots,j_{n}=1}^{k_1,\ldots,k_n}w_{j_1}\ldots w_{j_n}\Phi_{j_1,\ldots,j_n}(\bm{A}_{j_1}),\ldots,\sum\limits_{j_{1}=1,\ldots,j_{n}=1}^{k_1,\ldots,k_n}w_{j_1}\ldots w_{j_n}\Phi_{j_1,\ldots,j_n}(\bm{A}_{j_n})\Bigg)
\Bigg)\nonumber \\
\leq_2 \max\limits_{m_i \leq x_i \leq M_i}F\Big(\sum\limits_{i=1}^n c_i x_i + d, g(x_1,x_2,\ldots,x_n)\Big)\bm{I}_{\mathfrak{K}},~~~~~~~~~~~~~~~~~~~~~~~~~~~~~~~~~~~~~~~~~~~~~~~~~~~~~~~~~
\end{eqnarray}
where the inequality $\leq_1$ comes from the operator monotonicity of $F(u,v)$ with respect to the variable $u$, and the inequality $\leq_2$ comes from Eq.~\eqref{eq4 UB: thm:main 2.3}. Then, we have the desired inequality provided by Eq.~\eqref{eq UB: thm:main 2.3}.

Now, we will prove Eq.~\eqref{eq LB: thm:main 2.3}. Since $f(x_1,x_2,\ldots,x_n) \geq \sum\limits_{i=1}^n a_i x_i + b$ for $[x_1,x_2,\ldots,x_n] \in \bigtimes\limits_{i=1}^n [m_i,M_i]$ from Eq.~\eqref{eq f bounds: thm:main 2.3}, we have
\begin{eqnarray}\label{eq1 LB: thm:main 2.3}
f(\bm{A}_{j_1},\bm{A}_{j_2},\ldots,\bm{A}_{j_n}) \geq \sum\limits_{i=1}^n a_i \bm{A}_{j_i} + b\bm{I}_{\mathfrak{H}},
\end{eqnarray}
where $j_i=1,2,\ldots,k_i$ for $i=1,2,\ldots,n$. Since $\Phi_{j_1,\ldots,j_n}$ is a normalized positive linear map,  we also have
\begin{eqnarray}\label{eq2 LB: thm:main 2.3}
\Phi_{j_1,\ldots,j_n}(f(\bm{A}_{j_1},\bm{A}_{j_2},\ldots,\bm{A}_{j_n}))&\geq&\Phi_{j_1,\ldots,j_n}\Big(\sum\limits_{i=1}^n a_i \bm{A}_{j_i} + b\bm{I}_{\mathfrak{H}}\Big)\nonumber \\
&=& \sum\limits_{i=1}^n a_i \Phi_{j_1,\ldots,j_n}(\bm{A}_{j_i}) + b\bm{I}_{\mathfrak{K}},
\end{eqnarray}
where $j_i=1,2,\ldots,k_i$ for $i=1,2,\ldots,n$. By multiplying $w_{j_1}\ldots w_{j_n}$ to both sides of Eq.~\eqref{eq2 LB: thm:main 2.3} and summing items with respect to $w_{j_1}\ldots w_{j_n}$, we have 
\begin{eqnarray}\label{eq3 LB: thm:main 2.3}
\lefteqn{\sum\limits_{j_{1}=1,\ldots,j_{n}=1}^{k_1,\ldots,k_n}  w_{j_1}\ldots w_{j_n}\Phi_{j_1,\ldots,j_n}(f(\bm{A}_{j_1},\bm{A}_{j_2},\ldots,\bm{A}_{j_n}))}\nonumber \\
&\geq&
\sum\limits_{i=1}^n c_i \sum\limits_{j_{1}=1,\ldots,j_{n}=1}^{k_1,\ldots,k_n}w_{j_1}\ldots w_{j_n}\Phi_{j_1,\ldots,j_n}(\bm{A}_{j_i})+d\bm{I}_{\mathfrak{K}}.
\end{eqnarray}
Because $m_i \bm{I}_{\mathfrak{H}} \leq \bm{A}_{j_i} \leq M_i \bm{I}_{\mathfrak{H}}$ where $j_i=1,2,\ldots,k_i$ for $i=1,2,\ldots,n$, we have 
\begin{eqnarray}\label{eq4 LB: thm:main 2.3}
\Lambda\left(\sum\limits_{j_{1}=1,\ldots,j_{n}=1}^{k_1,\ldots,k_n}w_{j_1}\ldots w_{j_n}\Phi_{j_1,\ldots,j_n}(\bm{A}_{j_i})\right) \in [m_i, M_i].
\end{eqnarray}
By utilizing the operator monotonicity of $F(u,v)$ with respect to the variable $u$, we have
\begin{eqnarray}\label{eq5 LB: thm:main 2.3}
F\Bigg(\sum\limits_{j_{1}=1,\ldots,j_{n}=1}^{k_1,\ldots,k_n}  w_{j_1}\ldots w_{j_n}\Phi_{j_1,\ldots,j_n}(f(\bm{A}_{j_1},\bm{A}_{j_2},\ldots,\bm{A}_{j_n})),~~~~~~~~~~~~~~~~~~~~~~~~~~~~~~~~~~~~~~~~~~~~~~~~~~\nonumber \\
g\Bigg(\sum\limits_{j_{1}=1,\ldots,j_{n}=1}^{k_1,\ldots,k_n}w_{j_1}\ldots w_{j_n}\Phi_{j_1,\ldots,j_n}(\bm{A}_{j_1}),\ldots,\sum\limits_{j_{1}=1,\ldots,j_{n}=1}^{k_1,\ldots,k_n}w_{j_1}\ldots w_{j_n}\Phi_{j_1,\ldots,j_n}(\bm{A}_{j_n})\Bigg)
\Bigg)\nonumber \\
\geq_1 F\Bigg(\sum\limits_{i=1}^n c_i \sum\limits_{j_{1}=1,\ldots,j_{n}=1}^{k_1,\ldots,k_n}w_{j_1}\ldots w_{j_n}\Phi_{j_1,\ldots,j_n}(\bm{A}_{j_i})+d\bm{I}_{\mathfrak{K}},~~~~~~~~~~~~~~~~~~~~~~~~~~~~~~~~~~~~~~~~~~~~~~~~~~\nonumber \\
g\Bigg(\sum\limits_{j_{1}=1,\ldots,j_{n}=1}^{k_1,\ldots,k_n}w_{j_1}\ldots w_{j_n}\Phi_{j_1,\ldots,j_n}(\bm{A}_{j_1}),\ldots,\sum\limits_{j_{1}=1,\ldots,j_{n}=1}^{k_1,\ldots,k_n}w_{j_1}\ldots w_{j_n}\Phi_{j_1,\ldots,j_n}(\bm{A}_{j_n})\Bigg)
\Bigg)\nonumber \\
\geq_2 \min\limits_{m_i \leq x_i \leq M_i}F\Big(\sum\limits_{i=1}^n c_i x_i + d, g(x_1,x_2,\ldots,x_n)\Big)\bm{I}_{\mathfrak{K}},~~~~~~~~~~~~~~~~~~~~~~~~~~~~~~~~~~~~~~~~~~~~~~~~~~~~~~~~~
\end{eqnarray}
where the inequality $\geq_1$ comes from the operator monotonicity of $F(u,v)$ with respect to the variable $u$, and the inequality $\geq_2$ comes from Eq.~\eqref{eq4 LB: thm:main 2.3}. Then, we have the desired inequality provided by Eq.~\eqref{eq LB: thm:main 2.3}. 
$\hfill \Box$

By setting the function $F(u,v)$ in the following format:
\begin{eqnarray}\label{eq: F u v }
F(u,v) &=& u - \alpha v,
\end{eqnarray} 
where $\alpha \in \mathbb{R}$, we establish Theorem~\ref{thm:main 2.4}. This theorem lays out multivariate hypercomplex functions inequalities with both ratio and difference variations, which will be elaborated upon in subsequent sections.

\begin{theorem}\label{thm:main 2.4}
Let $\bm{A}_{j_i}$ be self-adjoint operators with $\Lambda(\bm{A}_{j_i}) \in [m_i, M_i]$ for real scalars $m_i <  M_i$. The mappings $\Phi_{j_1,\ldots,j_n}: \mathscr{B}(\mathfrak{H}) \rightarrow \mathscr{B}(\mathfrak{K})$ are normalized positive linear maps, where $j_i=1,2,\ldots,k_i$ for $i=1,2,\ldots,n$. We have $n$ probability vectors $\bm{w}_i =[w_{i,1},w_{i,2},\cdots, w_{i,k_i}]$ with the dimension $k_i$ for $i=1,2,\ldots,n$, i.e., $\sum\limits_{\ell=1}^{k_i}w_{i,\ell} = 1$. Let $f(x_1,x_2,\ldots,x_n)$ be any real valued continuous functions with $n$ variables defined on the range $\bigtimes\limits_{i=1}^n [m_i, M_i] \in \mathbb{R}^n$, where $\times$ is the Cartesian product. Besides, we assume that the function $f$ satisfies the following:
\begin{eqnarray}\label{eq f bounds: thm:main 2.4}
\sum\limits_{i=1}^n a_i x_i + b \leq f(x_1,x_2,\ldots,x_n) \leq \sum\limits_{i=1}^n c_i x_i + d,
\end{eqnarray}
for $[x_1,x_2,\ldots,x_n] \in \bigtimes\limits_{i=1}^n [m_i,M_i]$. The function $g(x_1,x_2,\ldots,x_n)$ is also a real valued continuous function with $n$ variables defined on the range $\bigtimes\limits_{i=1}^n [m_i, M_i]$. We also have a real valued function $F(u,v)$ defined by Eq.~\eqref{eq: F u v } with operator monotone on the first variable $u$ defined on $U \times V$ such that $f(\bigtimes\limits_{i=1}^n [m_i, M_i]) \subset U$, and $g(\bigtimes\limits_{i=1}^n [m_i, M_i]) \subset V$. Then, we have the following upper bound:
\begin{eqnarray}\label{eq UB: thm:main 2.4}
\lefteqn{\sum\limits_{j_{1}=1,\ldots,j_{n}=1}^{k_1,\ldots,k_n}  w_{j_1}\ldots w_{j_n}\Phi_{j_1,\ldots,j_n}(f(\bm{A}_{j_1},\bm{A}_{j_2},\ldots,\bm{A}_{j_n}))}\nonumber \\
&\leq&\alpha g\Bigg(\sum\limits_{j_{1}=1,\ldots,j_{n}=1}^{k_1,\ldots,k_n}w_{j_1}\ldots w_{j_n}\Phi_{j_1,\ldots,j_n}(\bm{A}_{j_1}),\ldots,\sum\limits_{j_{1}=1,\ldots,j_{n}=1}^{k_1,\ldots,k_n}w_{j_1}\ldots w_{j_n}\Phi_{j_1,\ldots,j_n}(\bm{A}_{j_n})\Bigg)
\nonumber \\
&&+
\max\limits_{\substack{m_i \leq x_i \leq M_i \\ i=1,2,\ldots,n}}\Big(\sum\limits_{i=1}^n c_i x_i + d - \alpha g(x_1,x_2,\ldots,x_n)\Big)\bm{I}_{\mathfrak{K}}.
\end{eqnarray}
Similarly, we also have the following lower bound:
\begin{eqnarray}\label{eq LB: thm:main 2.4}
\lefteqn{\sum\limits_{j_{1}=1,\ldots,j_{n}=1}^{k_1,\ldots,k_n}  w_{j_1}\ldots w_{j_n}\Phi_{j_1,\ldots,j_n}(f(\bm{A}_{j_1},\bm{A}_{j_2},\ldots,\bm{A}_{j_n}))}\nonumber \\
&\geq&\alpha g\Bigg(\sum\limits_{j_{1}=1,\ldots,j_{n}=1}^{k_1,\ldots,k_n}w_{j_1}\ldots w_{j_n}\Phi_{j_1,\ldots,j_n}(\bm{A}_{j_1}),\ldots,\sum\limits_{j_{1}=1,\ldots,j_{n}=1}^{k_1,\ldots,k_n}w_{j_1}\ldots w_{j_n}\Phi_{j_1,\ldots,j_n}(\bm{A}_{j_n})\Bigg)
\nonumber \\
&&+
\min\limits_{\substack{m_i \leq x_i \leq M_i \\ i=1,2,\ldots,n}}\Big(\sum\limits_{i=1}^n a_i x_i + b - \alpha g(x_1,x_2,\ldots,x_n)\Big)\bm{I}_{\mathfrak{K}}.
\end{eqnarray}
\end{theorem}
\textbf{Proof:}
By setting $F(u,v)=u - \alpha v$ in Eq.~\eqref{eq UB: thm:main 2.3} in Theorem~\ref{thm: main 2.3}, we have the desired inequality provided by Eq.~\eqref{eq UB: thm:main 2.4}. Similarly, By setting $F(u,v)=u - \alpha v$ in Eq.~\eqref{eq LB: thm:main 2.3} in Theorem~\ref{thm: main 2.3}, we have the desired inequality provided by Eq.~\eqref{eq LB: thm:main 2.4}.
$\hfill \Box$

Corollary~\ref{cor: special cases of g func} below is provided to give upper and lower bounds for special types of the function $g$ by applying Theorem~\ref{thm:main 2.4}.

\begin{corollary}\label{cor: special cases of g func}
Let $\bm{A}_{j_i}$ be self-adjoint operators with $\Lambda(\bm{A}_{j_i}) \in [m_i, M_i]$ for real scalars $m_i <  M_i$. The mappings $\Phi_{j_1,\ldots,j_n}: \mathscr{B}(\mathfrak{H}) \rightarrow \mathscr{B}(\mathfrak{K})$ are normalized positive linear maps, where $j_i=1,2,\ldots,k_i$ for $i=1,2,\ldots,n$. We have $n$ probability vectors $\bm{w}_i =[w_{i,1},w_{i,2},\cdots, w_{i,k_i}]$ with the dimension $k_i$ for $i=1,2,\ldots,n$, i.e., $\sum\limits_{\ell=1}^{k_i}w_{i,\ell} = 1$. Let $f(x_1,x_2,\ldots,x_n)$ be any real valued continuous functions with $n$ variables defined on the range $\bigtimes\limits_{i=1}^n [m_i, M_i] \in \mathbb{R}^n$, where $\times$ is the Cartesian product. Besides, we assume that the function $f$ satisfies the following:
\begin{eqnarray}\label{eq f bounds: cor: special cases of g func}
\sum\limits_{i=1}^n a_i x_i + b \leq f(x_1,x_2,\ldots,x_n) \leq \sum\limits_{i=1}^n c_i x_i + d,
\end{eqnarray}
for $[x_1,x_2,\ldots,x_n] \in \bigtimes\limits_{i=1}^n [m_i,M_i]$. 

(I) If $g(x_1,\ldots,x_n) = \Big(\sum\limits_{i=1}^n \beta_i x_i\Big)^q$, where $\beta_i \geq 0$, $q \in \mathbb{R}$ and $m_i \geq 0$, we have the upper bound for $\sum\limits_{j_{1}=1,\ldots,j_{n}=1}^{k_1,\ldots,k_n}  w_{j_1}\ldots w_{j_n}\Phi_{j_1,\ldots,j_n}(f(\bm{A}_{j_1},\bm{A}_{j_2},\ldots,\bm{A}_{j_n}))$:
\begin{eqnarray}\label{eq power U: cor: special cases of g func}
\lefteqn{\sum\limits_{j_{1}=1,\ldots,j_{n}=1}^{k_1,\ldots,k_n}  w_{j_1}\ldots w_{j_n}\Phi_{j_1,\ldots,j_n}(f(\bm{A}_{j_1},\bm{A}_{j_2},\ldots,\bm{A}_{j_n}))}\nonumber \\
&\leq&\alpha \left(\sum\limits_{i=1}^n \beta_i  \sum\limits_{j_{1}=1,\ldots,j_{n}=1}^{k_1,\ldots,k_n}w_{j_1}\ldots w_{j_n}\Phi_{j_1,\ldots,j_n}(\bm{A}_{j_i})\right)^q \nonumber \\
& &+
\max\limits_{\substack{m_i \leq x_i \leq M_i \\ i=1,2,\ldots,n}}\Big[\sum\limits_{i=1}^n c_i x_i + d - \alpha \Big(\sum\limits_{i=1}^n \beta_i x_i\Big)^q\Big]\bm{I}_{\mathfrak{K}},
\end{eqnarray}
and we have the lower bound for $\sum\limits_{j_{1}=1,\ldots,j_{n}=1}^{k_1,\ldots,k_n}  w_{j_1}\ldots w_{j_n}\Phi_{j_1,\ldots,j_n}(f(\bm{A}_{j_1},\bm{A}_{j_2},\ldots,\bm{A}_{j_n}))$:
\begin{eqnarray}\label{eq power L: cor: special cases of g func}
\lefteqn{\sum\limits_{j_{1}=1,\ldots,j_{n}=1}^{k_1,\ldots,k_n}  w_{j_1}\ldots w_{j_n}\Phi_{j_1,\ldots,j_n}(f(\bm{A}_{j_1},\bm{A}_{j_2},\ldots,\bm{A}_{j_n}))}\nonumber \\
&\geq&\alpha \left(\sum\limits_{i=1}^n \beta_i  \sum\limits_{j_{1}=1,\ldots,j_{n}=1}^{k_1,\ldots,k_n}w_{j_1}\ldots w_{j_n}\Phi_{j_1,\ldots,j_n}(\bm{A}_{j_i})\right)^q \nonumber \\
& &+
\min\limits_{\substack{m_i \leq x_i \leq M_i \\ i=1,2,\ldots,n}}\Big[\sum\limits_{i=1}^n a_i x_i + b - \alpha \Big(\sum\limits_{i=1}^n \beta_i x_i\Big)^q\Big]\bm{I}_{\mathfrak{K}}.
\end{eqnarray}

(II) If $g(x_1,\ldots,x_n) = \log\Big(\sum\limits_{i=1}^n \beta_i x_i\Big)$, where $\beta_i \geq 0$, and $m_i > 0$, we have the upper bound for $\sum\limits_{j_{1}=1,\ldots,j_{n}=1}^{k_1,\ldots,k_n}  w_{j_1}\ldots w_{j_n}\Phi_{j_1,\ldots,j_n}(f(\bm{A}_{j_1},\bm{A}_{j_2},\ldots,\bm{A}_{j_n}))$:
\begin{eqnarray}\label{eq log U: cor: special cases of g func}
\lefteqn{\sum\limits_{j_{1}=1,\ldots,j_{n}=1}^{k_1,\ldots,k_n}  w_{j_1}\ldots w_{j_n}\Phi_{j_1,\ldots,j_n}(f(\bm{A}_{j_1},\bm{A}_{j_2},\ldots,\bm{A}_{j_n}))}\nonumber \\
&\leq&\alpha\log\left(\sum\limits_{i=1}^n \beta_i  \sum\limits_{j_{1}=1,\ldots,j_{n}=1}^{k_1,\ldots,k_n}w_{j_1}\ldots w_{j_n}\Phi_{j_1,\ldots,j_n}(\bm{A}_{j_i})\right) \nonumber \\
& &+
\max\limits_{\substack{m_i \leq x_i \leq M_i \\ i=1,2,\ldots,n}}\Big[\sum\limits_{i=1}^n c_i x_i + d - \alpha\log\Big(\sum\limits_{i=1}^n \beta_i x_i\Big)\Big]\bm{I}_{\mathfrak{K}},
\end{eqnarray}
and we have the lower bound for $\sum\limits_{j_{1}=1,\ldots,j_{n}=1}^{k_1,\ldots,k_n}  w_{j_1}\ldots w_{j_n}\Phi_{j_1,\ldots,j_n}(f(\bm{A}_{j_1},\bm{A}_{j_2},\ldots,\bm{A}_{j_n}))$:
\begin{eqnarray}\label{eq log L: cor: special cases of g func}
\lefteqn{\sum\limits_{j_{1}=1,\ldots,j_{n}=1}^{k_1,\ldots,k_n}  w_{j_1}\ldots w_{j_n}\Phi_{j_1,\ldots,j_n}(f(\bm{A}_{j_1},\bm{A}_{j_2},\ldots,\bm{A}_{j_n}))}\nonumber \\
&\geq&\alpha\log\left(\sum\limits_{i=1}^n \beta_i  \sum\limits_{j_{1}=1,\ldots,j_{n}=1}^{k_1,\ldots,k_n}w_{j_1}\ldots w_{j_n}\Phi_{j_1,\ldots,j_n}(\bm{A}_{j_i})\right) \nonumber \\
& &+
\min\limits_{\substack{m_i \leq x_i \leq M_i \\ i=1,2,\ldots,n}}\Big[\sum\limits_{i=1}^n a_i x_i + b - \alpha\log\Big(\sum\limits_{i=1}^n \beta_i x_i\Big)\Big]\bm{I}_{\mathfrak{K}}.
\end{eqnarray}

(III) If $g(x_1,\ldots,x_n) = \exp\Big(\sum\limits_{i=1}^n \beta_i x_i\Big)$, we have the upper bound for \\
$\sum\limits_{j_{1}=1,\ldots,j_{n}=1}^{k_1,\ldots,k_n}  w_{j_1}\ldots w_{j_n}\Phi_{j_1,\ldots,j_n}(f(\bm{A}_{j_1},\bm{A}_{j_2},\ldots,\bm{A}_{j_n}))$:
\begin{eqnarray}\label{eq exp U: cor: special cases of g func}
\lefteqn{\sum\limits_{j_{1}=1,\ldots,j_{n}=1}^{k_1,\ldots,k_n}  w_{j_1}\ldots w_{j_n}\Phi_{j_1,\ldots,j_n}(f(\bm{A}_{j_1},\bm{A}_{j_2},\ldots,\bm{A}_{j_n}))}\nonumber \\
&\leq&\alpha\exp\left(\sum\limits_{i=1}^n \beta_i  \sum\limits_{j_{1}=1,\ldots,j_{n}=1}^{k_1,\ldots,k_n}w_{j_1}\ldots w_{j_n}\Phi_{j_1,\ldots,j_n}(\bm{A}_{j_i})\right) \nonumber \\
& &+
\max\limits_{\substack{m_i \leq x_i \leq M_i \\ i=1,2,\ldots,n}}\Big[\sum\limits_{i=1}^n c_i x_i + d - \alpha\exp\Big(\sum\limits_{i=1}^n \beta_i x_i\Big)\Big]\bm{I}_{\mathfrak{K}},
\end{eqnarray}
and we have the lower bound for $\sum\limits_{j_{1}=1,\ldots,j_{n}=1}^{k_1,\ldots,k_n}  w_{j_1}\ldots w_{j_n}\Phi_{j_1,\ldots,j_n}(f(\bm{A}_{j_1},\bm{A}_{j_2},\ldots,\bm{A}_{j_n}))$:
\begin{eqnarray}\label{eq exp L: cor: special cases of g func}
\lefteqn{\sum\limits_{j_{1}=1,\ldots,j_{n}=1}^{k_1,\ldots,k_n}  w_{j_1}\ldots w_{j_n}\Phi_{j_1,\ldots,j_n}(f(\bm{A}_{j_1},\bm{A}_{j_2},\ldots,\bm{A}_{j_n}))}\nonumber \\
&\geq&\alpha\exp\left(\sum\limits_{i=1}^n \beta_i  \sum\limits_{j_{1}=1,\ldots,j_{n}=1}^{k_1,\ldots,k_n}w_{j_1}\ldots w_{j_n}\Phi_{j_1,\ldots,j_n}(\bm{A}_{j_i})\right) \nonumber \\
& &+
\min\limits_{\substack{m_i \leq x_i \leq M_i \\ i=1,2,\ldots,n}}\Big[\sum\limits_{i=1}^n a_i x_i + b - \alpha\exp\Big(\sum\limits_{i=1}^n \beta_i x_i\Big)\Big]\bm{I}_{\mathfrak{K}}.
\end{eqnarray}
\end{corollary}
\textbf{Proof:}
This corollary is obtained by applying Theprem~\ref{thm:main 2.4} directly to functions $g$ as 
$g(x_1,\ldots,x_n) = \Big(\sum\limits_{i=1}^n \beta_i x_i\Big)^q$, $g(x_1,\ldots,x_n) = \log\Big(\sum\limits_{i=1}^n \beta_i x_i\Big)$ and $g(x_1,\ldots,x_n) = \exp\Big(\sum\limits_{i=1}^n \beta_i x_i\Big)$, respectivey.
$\hfill \Box$

\section{Multivariate Hypercomplex Functions Inequalities: Ratio Kind}\label{sec: Multivariate Hypercomplex Functions Inequalities: Ratio Kind}

This section will explore the determination of lower and upper bounds for \\
$\sum\limits_{j_{1}=1,\ldots,j_{n}=1}^{k_1,\ldots,k_n}  w_{j_1}\ldots w_{j_n}\Phi_{j_1,\ldots,j_n}(f(\bm{A}_{j_1},\bm{A}_{j_2},\ldots,\bm{A}_{j_n}))$ based on ratio criteria associated with the function $g(x_1, x_2,\ldots,x_n)$.

\begin{theorem}\label{thm: 2.9}
Let $\bm{A}_{j_i}$ be self-adjoint operators with $\Lambda(\bm{A}_{j_i}) \in [m_i, M_i]$ for real scalars $m_i <  M_i$. The mappings $\Phi_{j_1,\ldots,j_n}: \mathscr{B}(\mathfrak{H}) \rightarrow \mathscr{B}(\mathfrak{K})$ are normalized positive linear maps, where $j_i=1,2,\ldots,k_i$ for $i=1,2,\ldots,n$. We have $n$ probability vectors $\bm{w}_i =[w_{i,1},w_{i,2},\cdots, w_{i,k_i}]$ with the dimension $k_i$ for $i=1,2,\ldots,n$, i.e., $\sum\limits_{\ell=1}^{k_i}w_{i,\ell} = 1$. Let $f(x_1,x_2,\ldots,x_n)$ be any real valued continuous functions with $n$ variables defined on the range $\bigtimes\limits_{i=1}^n [m_i, M_i] \in \mathbb{R}^n$, where $\times$ is the Cartesian product. Besides, we assume that the function $f$ satisfies the following:
\begin{eqnarray}\label{eq f bounds: thm: 2.9}
\sum\limits_{i=1}^n a_i x_i + b \leq f(x_1,x_2,\ldots,x_n) \leq \sum\limits_{i=1}^n c_i x_i + d,
\end{eqnarray}
for $[x_1,x_2,\ldots,x_n] \in \bigtimes\limits_{i=1}^n [m_i,M_i]$.

(I) If we also assume that $g(x_1,\ldots,x_n) > 0$ for $[x_1,\ldots,x_n] \in \bigtimes\limits_{i=1}^n [m_i, M_i]$, then, we have the following upper bound for $\sum\limits_{j_{1}=1,\ldots,j_{n}=1}^{k_1,\ldots,k_n}  w_{j_1}\ldots w_{j_n}\Phi_{j_1,\ldots,j_n}(f(\bm{A}_{j_1},\bm{A}_{j_2},\ldots,\bm{A}_{j_n}))$:
\begin{eqnarray}\label{eq1 UB: thm: 2.9 pos g}
\lefteqn{\sum\limits_{j_{1}=1,\ldots,j_{n}=1}^{k_1,\ldots,k_n}  w_{j_1}\ldots w_{j_n}\Phi_{j_1,\ldots,j_n}(f(\bm{A}_{j_1},\bm{A}_{j_2},\ldots,\bm{A}_{j_n}))}\nonumber \\
&\leq&\left[\max\limits_{\substack{m_i \leq x_i \leq M_i \\ i=1,2,\ldots,n}}\frac{\sum\limits_{i=1}^n c_i x_i +d}{g(x_1,\ldots,x_n)}\right]\nonumber \\
&&\times g\Bigg(\sum\limits_{j_{1}=1,\ldots,j_{n}=1}^{k_1,\ldots,k_n}w_{j_1}\ldots w_{j_n}\Phi_{j_1,\ldots,j_n}(\bm{A}_{j_1}),\ldots,\sum\limits_{j_{1}=1,\ldots,j_{n}=1}^{k_1,\ldots,k_n}w_{j_1}\ldots w_{j_n}\Phi_{j_1,\ldots,j_n}(\bm{A}_{j_n})\Bigg);
\end{eqnarray}
and, the following lower bound for $\sum\limits_{j_{1}=1,\ldots,j_{n}=1}^{k_1,\ldots,k_n}  w_{j_1}\ldots w_{j_n}\Phi_{j_1,\ldots,j_n}(f(\bm{A}_{j_1},\bm{A}_{j_2},\ldots,\bm{A}_{j_n}))$:
\begin{eqnarray}\label{eq1 LB: thm: 2.9 pos g}
\lefteqn{\sum\limits_{j_{1}=1,\ldots,j_{n}=1}^{k_1,\ldots,k_n}  w_{j_1}\ldots w_{j_n}\Phi_{j_1,\ldots,j_n}(f(\bm{A}_{j_1},\bm{A}_{j_2},\ldots,\bm{A}_{j_n}))}\nonumber \\
&\geq&\left[\min\limits_{\substack{m_i \leq x_i \leq M_i \\ i=1,2,\ldots,n}}\frac{\sum\limits_{i=1}^n a_i x_i +b}{g(x_1,\ldots,x_n)}\right]\nonumber \\
&&\times g\Bigg(\sum\limits_{j_{1}=1,\ldots,j_{n}=1}^{k_1,\ldots,k_n}w_{j_1}\ldots w_{j_n}\Phi_{j_1,\ldots,j_n}(\bm{A}_{j_1}),\ldots,\sum\limits_{j_{1}=1,\ldots,j_{n}=1}^{k_1,\ldots,k_n}w_{j_1}\ldots w_{j_n}\Phi_{j_1,\ldots,j_n}(\bm{A}_{j_n})\Bigg).
\end{eqnarray}

(II) If we also assume that $g(x_1,\ldots,x_n) < 0$ for $[x_1,\ldots,x_n] \in \bigtimes\limits_{i=1}^n [m_i, M_i]$, then, we have the following upper bound for $\sum\limits_{j_{1}=1,\ldots,j_{n}=1}^{k_1,\ldots,k_n}  w_{j_1}\ldots w_{j_n}\Phi_{j_1,\ldots,j_n}(f(\bm{A}_{j_1},\bm{A}_{j_2},\ldots,\bm{A}_{j_n}))$:
\begin{eqnarray}\label{eq1 UB: thm: 2.9 neg g}
\lefteqn{\sum\limits_{j_{1}=1,\ldots,j_{n}=1}^{k_1,\ldots,k_n}  w_{j_1}\ldots w_{j_n}\Phi_{j_1,\ldots,j_n}(f(\bm{A}_{j_1},\bm{A}_{j_2},\ldots,\bm{A}_{j_n}))}\nonumber \\
&\leq&\left[\min\limits_{\substack{m_i \leq x_i \leq M_i \\ i=1,2,\ldots,n}}\frac{\sum\limits_{i=1}^n c_i x_i +d}{g(x_1,\ldots,x_n)}\right]\nonumber \\
&&\times g\Bigg(\sum\limits_{j_{1}=1,\ldots,j_{n}=1}^{k_1,\ldots,k_n}w_{j_1}\ldots w_{j_n}\Phi_{j_1,\ldots,j_n}(\bm{A}_{j_1}),\ldots,\sum\limits_{j_{1}=1,\ldots,j_{n}=1}^{k_1,\ldots,k_n}w_{j_1}\ldots w_{j_n}\Phi_{j_1,\ldots,j_n}(\bm{A}_{j_n})\Bigg);
\end{eqnarray}
and, the following lower bound for $\sum\limits_{j_{1}=1,\ldots,j_{n}=1}^{k_1,\ldots,k_n}  w_{j_1}\ldots w_{j_n}\Phi_{j_1,\ldots,j_n}(f(\bm{A}_{j_1},\bm{A}_{j_2},\ldots,\bm{A}_{j_n}))$:
\begin{eqnarray}\label{eq1 LB: thm: 2.9 neg g}
\lefteqn{\sum\limits_{j_{1}=1,\ldots,j_{n}=1}^{k_1,\ldots,k_n}  w_{j_1}\ldots w_{j_n}\Phi_{j_1,\ldots,j_n}(f(\bm{A}_{j_1},\bm{A}_{j_2},\ldots,\bm{A}_{j_n}))}\nonumber \\
&\geq&\left[\max\limits_{\substack{m_i \leq x_i \leq M_i \\ i=1,2,\ldots,n}}\frac{\sum\limits_{i=1}^n a_i x_i +b}{g(x_1,\ldots,x_n)}\right]\nonumber \\
&&\times g\Bigg(\sum\limits_{j_{1}=1,\ldots,j_{n}=1}^{k_1,\ldots,k_n}w_{j_1}\ldots w_{j_n}\Phi_{j_1,\ldots,j_n}(\bm{A}_{j_1}),\ldots,\sum\limits_{j_{1}=1,\ldots,j_{n}=1}^{k_1,\ldots,k_n}w_{j_1}\ldots w_{j_n}\Phi_{j_1,\ldots,j_n}(\bm{A}_{j_n})\Bigg).
\end{eqnarray}
\end{theorem}
\textbf{Proof:}
For part (I), we will apply $F(u,v)$ as 
\begin{eqnarray}\label{eq: F u v 2}
F(u,v)&=&v^{-1/2} u v^{-1/2},
\end{eqnarray}
to Eq.~\eqref{eq UB: thm:main 2.3} in Theorem~\ref{thm: main 2.3}, then, we will obtain 
\begin{eqnarray}\label{eq2 UB: thm: 2.9 pos g}
\lefteqn{g\Bigg(\sum\limits_{j_{1}=1,\ldots,j_{n}=1}^{k_1,\ldots,k_n}w_{j_1}\ldots w_{j_n}\Phi_{j_1,\ldots,j_n}(\bm{A}_{j_1}),\ldots,\sum\limits_{j_{1}=1,\ldots,j_{n}=1}^{k_1,\ldots,k_n}w_{j_1}\ldots w_{j_n}\Phi_{j_1,\ldots,j_n}(\bm{A}_{j_n})\Bigg)^{-1/2}}\nonumber \\
&& \times
\left(\sum\limits_{j_{1}=1,\ldots,j_{n}=1}^{k_1,\ldots,k_n}  w_{j_1}\ldots w_{j_n}\Phi_{j_1,\ldots,j_n}(f(\bm{A}_{j_1},\bm{A}_{j_2},\ldots,\bm{A}_{j_n}))\right) \nonumber \\
&& \times g\Bigg(\sum\limits_{j_{1}=1,\ldots,j_{n}=1}^{k_1,\ldots,k_n}w_{j_1}\ldots w_{j_n}\Phi_{j_1,\ldots,j_n}(\bm{A}_{j_1}),\ldots,\sum\limits_{j_{1}=1,\ldots,j_{n}=1}^{k_1,\ldots,k_n}w_{j_1}\ldots w_{j_n}\Phi_{j_1,\ldots,j_n}(\bm{A}_{j_n})\Bigg)^{-1/2}
\nonumber \\
&\leq&\left[\max\limits_{\substack{m_i \leq x_i \leq M_i \\ i=1,2,\ldots,n}}\frac{\sum\limits_{i=1}^n c_i x_i +d}{g(x_1,\ldots,x_n)}\right]\bm{1}_{\mathfrak{K}}.
\end{eqnarray}
By multiplying $g\Bigg(\sum\limits_{j_{1}=1,\ldots,j_{n}=1}^{k_1,\ldots,k_n}w_{j_1}\ldots w_{j_n}\Phi_{j_1,\ldots,j_n}(\bm{A}_{j_1}),\ldots,\sum\limits_{j_{1}=1,\ldots,j_{n}=1}^{k_1,\ldots,k_n}w_{j_1}\ldots w_{j_n}\Phi_{j_1,\ldots,j_n}(\bm{A}_{j_n})\Bigg)^{1/2}$ at both sides of Eq.~\eqref{eq2 UB: thm: 2.9 pos g}, we obtain the desired inequality provided by Eq.~\eqref{eq1 UB: thm: 2.9 pos g}. By applying $F(u,v)$ with Eq.~\eqref{eq: F u v 2} again to Eq.~\eqref{eq LB: thm:main 2.3} in Theorem~\ref{thm: main 2.3}, then, we will obtain 
\begin{eqnarray}\label{eq2 LB: thm: 2.9 pos g}
\lefteqn{g\Bigg(\sum\limits_{j_{1}=1,\ldots,j_{n}=1}^{k_1,\ldots,k_n}w_{j_1}\ldots w_{j_n}\Phi_{j_1,\ldots,j_n}(\bm{A}_{j_1}),\ldots,\sum\limits_{j_{1}=1,\ldots,j_{n}=1}^{k_1,\ldots,k_n}w_{j_1}\ldots w_{j_n}\Phi_{j_1,\ldots,j_n}(\bm{A}_{j_n})\Bigg)^{-1/2}}\nonumber \\
&& \times
\left(\sum\limits_{j_{1}=1,\ldots,j_{n}=1}^{k_1,\ldots,k_n}  w_{j_1}\ldots w_{j_n}\Phi_{j_1,\ldots,j_n}(f(\bm{A}_{j_1},\bm{A}_{j_2},\ldots,\bm{A}_{j_n}))\right) \nonumber \\
&& \times g\Bigg(\sum\limits_{j_{1}=1,\ldots,j_{n}=1}^{k_1,\ldots,k_n}w_{j_1}\ldots w_{j_n}\Phi_{j_1,\ldots,j_n}(\bm{A}_{j_1}),\ldots,\sum\limits_{j_{1}=1,\ldots,j_{n}=1}^{k_1,\ldots,k_n}w_{j_1}\ldots w_{j_n}\Phi_{j_1,\ldots,j_n}(\bm{A}_{j_n})\Bigg)^{-1/2}
\nonumber \\
&\geq&\left[\min\limits_{\substack{m_i \leq x_i \leq M_i \\ i=1,2,\ldots,n}}\frac{\sum\limits_{i=1}^n a_i x_i +b}{g(x_1,\ldots,x_n)}\right]\bm{1}_{\mathfrak{K}}.
\end{eqnarray}
By multiplying $g\Bigg(\sum\limits_{j_{1}=1,\ldots,j_{n}=1}^{k_1,\ldots,k_n}w_{j_1}\ldots w_{j_n}\Phi_{j_1,\ldots,j_n}(\bm{A}_{j_1}),\ldots,\sum\limits_{j_{1}=1,\ldots,j_{n}=1}^{k_1,\ldots,k_n}w_{j_1}\ldots w_{j_n}\Phi_{j_1,\ldots,j_n}(\bm{A}_{j_n})\Bigg)^{1/2}$ at both sides of Eq.~\eqref{eq2 LB: thm: 2.9 pos g}, we obtain the desired inequality provided by Eq.~\eqref{eq1 LB: thm: 2.9 pos g}. 

The proof of Part (II) is immediate obtained by setting $g 
(x_1,\ldots,x_n)$ as $-g(x_1,\ldots,x_n)$ in Eq.~\eqref{eq2 UB: thm: 2.9 pos g} and Eq.~\eqref{eq2 LB: thm: 2.9 pos g}. 
$\hfill \Box$

Next Corollary~\ref{cor: 2.11-14} is obtained by applying Theorem~\ref{thm: 2.9} to special types of function $g$. 
\begin{corollary}\label{cor: 2.11-14}
Let $\bm{A}_{j_i}$ be self-adjoint operators with $\Lambda(\bm{A}_{j_i}) \in [m_i, M_i]$ for real scalars $m_i <  M_i$. The mappings $\Phi_{j_1,\ldots,j_n}: \mathscr{B}(\mathfrak{H}) \rightarrow \mathscr{B}(\mathfrak{K})$ are normalized positive linear maps, where $j_i=1,2,\ldots,k_i$ for $i=1,2,\ldots,n$. We have $n$ probability vectors $\bm{w}_i =[w_{i,1},w_{i,2},\cdots, w_{i,k_i}]$ with the dimension $k_i$ for $i=1,2,\ldots,n$, i.e., $\sum\limits_{\ell=1}^{k_i}w_{i,\ell} = 1$. Let $f(x_1,x_2,\ldots,x_n)$ be any real valued continuous functions with $n$ variables defined on the range $\bigtimes\limits_{i=1}^n [m_i, M_i] \in \mathbb{R}^n$, where $\times$ is the Cartesian product. Besides, we assume that the function $f$ satisfies the following:
\begin{eqnarray}\label{eq f bounds: cor: 2.11-14}
\sum\limits_{i=1}^n a_i x_i + b \leq f(x_1,x_2,\ldots,x_n) \leq \sum\limits_{i=1}^n c_i x_i + d,
\end{eqnarray}
for $[x_1,x_2,\ldots,x_n] \in \bigtimes\limits_{i=1}^n [m_i,M_i]$.

(I) If we also assume that $g(x_1,\ldots,x_n) = \Big(\sum\limits_{i=1}^n \beta_i x_i\Big)^q$, where $\beta_i \geq 0$, $q \in \mathbb{R}$ and $m_i \geq 0$, then, we have the following upper bound for $\sum\limits_{j_{1}=1,\ldots,j_{n}=1}^{k_1,\ldots,k_n}  w_{j_1}\ldots w_{j_n}\Phi_{j_1,\ldots,j_n}(f(\bm{A}_{j_1},\bm{A}_{j_2},\ldots,\bm{A}_{j_n}))$:
\begin{eqnarray}\label{eq x q UB: cor: 2.11-14}
\lefteqn{\sum\limits_{j_{1}=1,\ldots,j_{n}=1}^{k_1,\ldots,k_n}  w_{j_1}\ldots w_{j_n}\Phi_{j_1,\ldots,j_n}(f(\bm{A}_{j_1},\bm{A}_{j_2},\ldots,\bm{A}_{j_n}))}\nonumber \\
&\leq&\left[\max\limits_{\substack{m_i \leq x_i \leq M_i \\ i=1,2,\ldots,n}}\frac{\sum\limits_{i=1}^n c_i x_i +d}{\Big(\sum\limits_{i=1}^n \beta_i x_i\Big)^q}\right]\left(\sum\limits_{i=1}^n \beta_i \sum\limits_{j_{1}=1,\ldots,j_{n}=1}^{k_1,\ldots,k_n}w_{j_1}\ldots w_{j_n}\Phi_{j_1,\ldots,j_n}(\bm{A}_{j_i})\right)^q;
\end{eqnarray}
and, the following lower bound for $\sum\limits_{j_{1}=1,\ldots,j_{n}=1}^{k_1,\ldots,k_n}  w_{j_1}\ldots w_{j_n}\Phi_{j_1,\ldots,j_n}(f(\bm{A}_{j_1},\bm{A}_{j_2},\ldots,\bm{A}_{j_n}))$:
\begin{eqnarray}\label{eq x q LB: cor: 2.11-14}
\lefteqn{\sum\limits_{j_{1}=1,\ldots,j_{n}=1}^{k_1,\ldots,k_n}  w_{j_1}\ldots w_{j_n}\Phi_{j_1,\ldots,j_n}(f(\bm{A}_{j_1},\bm{A}_{j_2},\ldots,\bm{A}_{j_n}))}\nonumber \\
&\geq&\left[\min\limits_{\substack{m_i \leq x_i \leq M_i \\ i=1,2,\ldots,n}}\frac{\sum\limits_{i=1}^n a_i x_i +b}{\Big(\sum\limits_{i=1}^n \beta_i x_i\Big)^q}\right]\left(\sum\limits_{i=1}^n \beta_i \sum\limits_{j_{1}=1,\ldots,j_{n}=1}^{k_1,\ldots,k_n}w_{j_1}\ldots w_{j_n}\Phi_{j_1,\ldots,j_n}(\bm{A}_{j_i})\right)^q.
\end{eqnarray}

(II) If $g(x_1,\ldots,x_n) = \log\Big(\sum\limits_{i=1}^n \beta_i x_i\Big)$ with $\log\Big(\sum\limits_{i=1}^n \beta_i x_i\Big)>0$ for $[x_1,\ldots,x_n]\in\bigtimes\limits_{i=1}^n [m_i, M_i]$, we have the upper bound for $\sum\limits_{j_{1}=1,\ldots,j_{n}=1}^{k_1,\ldots,k_n}  w_{j_1}\ldots w_{j_n}\Phi_{j_1,\ldots,j_n}(f(\bm{A}_{j_1},\bm{A}_{j_2},\ldots,\bm{A}_{j_n}))$:
\begin{eqnarray}\label{eq log x l0 UB: cor: 2.11-14}
\lefteqn{\sum\limits_{j_{1}=1,\ldots,j_{n}=1}^{k_1,\ldots,k_n}  w_{j_1}\ldots w_{j_n}\Phi_{j_1,\ldots,j_n}(f(\bm{A}_{j_1},\bm{A}_{j_2},\ldots,\bm{A}_{j_n}))}\nonumber \\
&\leq&\left[\max\limits_{\substack{m_i \leq x_i \leq M_i \\ i=1,2,\ldots,n}}\frac{\sum\limits_{i=1}^n c_i x_i +d}{\log\Big(\sum\limits_{i=1}^n \beta_i x_i\Big)}\right]\log\left(\sum\limits_{i=1}^n \beta_i \sum\limits_{j_{1}=1,\ldots,j_{n}=1}^{k_1,\ldots,k_n}w_{j_1}\ldots w_{j_n}\Phi_{j_1,\ldots,j_n}(\bm{A}_{j_i})\right);
\end{eqnarray}
and, the following lower bound for $\sum\limits_{j=1}^k w_j\Phi(f(\bm{A}_j))$:
\begin{eqnarray}\label{eq log x l0 LB: cor: 2.11-14}
\lefteqn{\sum\limits_{j_{1}=1,\ldots,j_{n}=1}^{k_1,\ldots,k_n}  w_{j_1}\ldots w_{j_n}\Phi_{j_1,\ldots,j_n}(f(\bm{A}_{j_1},\bm{A}_{j_2},\ldots,\bm{A}_{j_n}))}\nonumber \\
&\geq&\left[\min\limits_{\substack{m_i \leq x_i \leq M_i \\ i=1,2,\ldots,n}}\frac{\sum\limits_{i=1}^n a_i x_i +b}{\log\Big(\sum\limits_{i=1}^n \beta_i x_i\Big)}\right]\log\left(\sum\limits_{i=1}^n \beta_i \sum\limits_{j_{1}=1,\ldots,j_{n}=1}^{k_1,\ldots,k_n}w_{j_1}\ldots w_{j_n}\Phi_{j_1,\ldots,j_n}(\bm{A}_{j_i})\right).
\end{eqnarray}

(II') If $g(x_1,\ldots,x_n) = \log\Big(\sum\limits_{i=1}^n \beta_i x_i\Big)$ with $\log\Big(\sum\limits_{i=1}^n \beta_i x_i\Big)<0$ for $[x_1,\ldots,x_n]\in\bigtimes\limits_{i=1}^n [m_i, M_i]$, we have the upper bound for $\sum\limits_{j_{1}=1,\ldots,j_{n}=1}^{k_1,\ldots,k_n}  w_{j_1}\ldots w_{j_n}\Phi_{j_1,\ldots,j_n}(f(\bm{A}_{j_1},\bm{A}_{j_2},\ldots,\bm{A}_{j_n}))$:
\begin{eqnarray}\label{eq log x s0 UB: cor: 2.11-14}
\lefteqn{\sum\limits_{j_{1}=1,\ldots,j_{n}=1}^{k_1,\ldots,k_n}  w_{j_1}\ldots w_{j_n}\Phi_{j_1,\ldots,j_n}(f(\bm{A}_{j_1},\bm{A}_{j_2},\ldots,\bm{A}_{j_n}))}\nonumber \\
&\leq&\left[\min\limits_{\substack{m_i \leq x_i \leq M_i \\ i=1,2,\ldots,n}}\frac{\sum\limits_{i=1}^n c_i x_i +d}{\log\Big(\sum\limits_{i=1}^n \beta_i x_i\Big)}\right]\log\left(\sum\limits_{i=1}^n \beta_i \sum\limits_{j_{1}=1,\ldots,j_{n}=1}^{k_1,\ldots,k_n}w_{j_1}\ldots w_{j_n}\Phi_{j_1,\ldots,j_n}(\bm{A}_{j_i})\right);
\end{eqnarray}
and, the following lower bound for $\sum\limits_{j=1}^k w_j\Phi(f(\bm{A}_j))$:
\begin{eqnarray}\label{eq log x s0 LB: cor: 2.11-14}
\lefteqn{\sum\limits_{j_{1}=1,\ldots,j_{n}=1}^{k_1,\ldots,k_n}  w_{j_1}\ldots w_{j_n}\Phi_{j_1,\ldots,j_n}(f(\bm{A}_{j_1},\bm{A}_{j_2},\ldots,\bm{A}_{j_n}))}\nonumber \\
&\geq&\left[\max\limits_{\substack{m_i \leq x_i \leq M_i \\ i=1,2,\ldots,n}}\frac{\sum\limits_{i=1}^n a_i x_i +b}{\log\Big(\sum\limits_{i=1}^n \beta_i x_i\Big)}\right]\log\left(\sum\limits_{i=1}^n \beta_i \sum\limits_{j_{1}=1,\ldots,j_{n}=1}^{k_1,\ldots,k_n}w_{j_1}\ldots w_{j_n}\Phi_{j_1,\ldots,j_n}(\bm{A}_{j_i})\right).
\end{eqnarray}

(III) If $g(x_1,\ldots,x_n) = \exp\Big(\sum\limits_{i=1}^n \beta_i x_i\Big)$, we have the upper bound for \\
$\sum\limits_{j_{1}=1,\ldots,j_{n}=1}^{k_1,\ldots,k_n}  w_{j_1}\ldots w_{j_n}\Phi_{j_1,\ldots,j_n}(f(\bm{A}_{j_1},\bm{A}_{j_2},\ldots,\bm{A}_{j_n}))$:
\begin{eqnarray}\label{eq exp x UB: cor: 2.11-14}
\lefteqn{\sum\limits_{j_{1}=1,\ldots,j_{n}=1}^{k_1,\ldots,k_n}  w_{j_1}\ldots w_{j_n}\Phi_{j_1,\ldots,j_n}(f(\bm{A}_{j_1},\bm{A}_{j_2},\ldots,\bm{A}_{j_n}))}\nonumber \\
&\leq&\left[\max\limits_{\substack{m_i \leq x_i \leq M_i \\ i=1,2,\ldots,n}}\frac{\sum\limits_{i=1}^n c_i x_i +d}{\exp\Big(\sum\limits_{i=1}^n \beta_i x_i\Big)}\right]\exp\left(\sum\limits_{i=1}^n \beta_i \sum\limits_{j_{1}=1,\ldots,j_{n}=1}^{k_1,\ldots,k_n}w_{j_1}\ldots w_{j_n}\Phi_{j_1,\ldots,j_n}(\bm{A}_{j_i})\right);
\end{eqnarray}
and, the following lower bound for  \\
$\sum\limits_{j_{1}=1,\ldots,j_{n}=1}^{k_1,\ldots,k_n}  w_{j_1}\ldots w_{j_n}\Phi_{j_1,\ldots,j_n}(f(\bm{A}_{j_1},\bm{A}_{j_2},\ldots,\bm{A}_{j_n}))$:
\begin{eqnarray}\label{eq exp x LB: cor: 2.11-14}
\lefteqn{\sum\limits_{j_{1}=1,\ldots,j_{n}=1}^{k_1,\ldots,k_n}  w_{j_1}\ldots w_{j_n}\Phi_{j_1,\ldots,j_n}(f(\bm{A}_{j_1},\bm{A}_{j_2},\ldots,\bm{A}_{j_n}))}\nonumber \\
&\geq&\left[\min\limits_{\substack{m_i \leq x_i \leq M_i \\ i=1,2,\ldots,n}}\frac{\sum\limits_{i=1}^n a_i x_i +b}{\exp\Big(\sum\limits_{i=1}^n \beta_i x_i\Big)}\right]\exp\left(\sum\limits_{i=1}^n \beta_i \sum\limits_{j_{1}=1,\ldots,j_{n}=1}^{k_1,\ldots,k_n}w_{j_1}\ldots w_{j_n}\Phi_{j_1,\ldots,j_n}(\bm{A}_{j_i})\right).
\end{eqnarray}
\end{corollary}
\textbf{Proof:}
Part (I) of this Corollary is proved by applying Theorem~\ref{thm: 2.9} Part (I) with the function $g$ as $g(x_1,\ldots,x_n) = \Big(\sum\limits_{i=1}^n \beta_i x_i\Big)^q$. Part (II) of this Corollary is proved by applying Theorem~\ref{thm: 2.9} Part (I) with the function $g$ as $g(x_1,\ldots,x_n) = \log\Big(\sum\limits_{i=1}^n \beta_i x_i\Big)$, where $\log\Big(\sum\limits_{i=1}^n \beta_i x_i\Big)> 0$. Part (II') of this Corollary is proved by applying Theorem~\ref{thm: 2.9} Part (II) with the function $g$ as $g(x_1,\ldots,x_n) = \log\Big(\sum\limits_{i=1}^n \beta_i x_i\Big)$, where $\log\Big(\sum\limits_{i=1}^n \beta_i x_i\Big)< 0$. Finally, Part (III) of this Corollary is proved by applying Theorem~\ref{thm: 2.9} Part (I) with the function $g$ as $g(x_1,\ldots,x_n) = \exp\Big(\sum\limits_{i=1}^n \beta_i x_i\Big)$.      
$\hfill \Box$

\section{Multivariate Hypercomplex Functions Inequalities: Difference Kind}\label{sec: Multivariate Hypercomplex Functions Inequalities: Difference Kind}

This section will focus on deriving the lower and upper bounds for $\sum\limits_{j_{1}=1,\ldots,j_{n}=1}^{k_1,\ldots,k_n}  w_{j_1}\ldots w_{j_n}\Phi_{j_1,\ldots,j_n}(f(\bm{A}_{j_1},\bm{A}_{j_2},\ldots,\bm{A}_{j_n}))$ using difference criteria associated with the function $g$.
\begin{theorem}\label{thm: cor 2.15}
Let $\bm{A}_{j_i}$ be self-adjoint operators with $\Lambda(\bm{A}_{j_i}) \in [m_i, M_i]$ for real scalars $m_i <  M_i$. The mappings $\Phi_{j_1,\ldots,j_n}: \mathscr{B}(\mathfrak{H}) \rightarrow \mathscr{B}(\mathfrak{K})$ are normalized positive linear maps, where $j_i=1,2,\ldots,k_i$ for $i=1,2,\ldots,n$. We have $n$ probability vectors $\bm{w}_i =[w_{i,1},w_{i,2},\cdots, w_{i,k_i}]$ with the dimension $k_i$ for $i=1,2,\ldots,n$, i.e., $\sum\limits_{\ell=1}^{k_i}w_{i,\ell} = 1$. Let $f(x_1,x_2,\ldots,x_n)$ be any real valued continuous functions with $n$ variables defined on the range $\bigtimes\limits_{i=1}^n [m_i, M_i] \in \mathbb{R}^n$, where $\times$ is the Cartesian product. Besides, we assume that the function $f$ satisfies the following:
\begin{eqnarray}\label{eq f bounds: thm: cor 2.15}
\sum\limits_{i=1}^n a_i x_i + b \leq f(x_1,x_2,\ldots,x_n) \leq \sum\limits_{i=1}^n c_i x_i + d,
\end{eqnarray}
for $[x_1,x_2,\ldots,x_n] \in \bigtimes\limits_{i=1}^n [m_i,M_i]$. 

Then, we have the following upper bound:
\begin{eqnarray}\label{eq UB: cor 2.15}
\lefteqn{\sum\limits_{j_{1}=1,\ldots,j_{n}=1}^{k_1,\ldots,k_n}  w_{j_1}\ldots w_{j_n}\Phi_{j_1,\ldots,j_n}(f(\bm{A}_{j_1},\bm{A}_{j_2},\ldots,\bm{A}_{j_n}))}\nonumber \\
&&- g\Bigg(\sum\limits_{j_{1}=1,\ldots,j_{n}=1}^{k_1,\ldots,k_n}w_{j_1}\ldots w_{j_n}\Phi_{j_1,\ldots,j_n}(\bm{A}_{j_1}),\ldots,\sum\limits_{j_{1}=1,\ldots,j_{n}=1}^{k_1,\ldots,k_n}w_{j_1}\ldots w_{j_n}\Phi_{j_1,\ldots,j_n}(\bm{A}_{j_n})\Bigg)
\nonumber \\
&\leq&
\max\limits_{\substack{m_i \leq x_i \leq M_i \\ i=1,2,\ldots,n}}\Big(\sum\limits_{i=1}^n c_i x_i + d - g(x_1,x_2,\ldots,x_n)\Big)\bm{I}_{\mathfrak{K}}.
\end{eqnarray}
Similarly, we also have the following lower bound:
\begin{eqnarray}\label{eq LB: cor 2.15}
\lefteqn{\sum\limits_{j_{1}=1,\ldots,j_{n}=1}^{k_1,\ldots,k_n}  w_{j_1}\ldots w_{j_n}\Phi_{j_1,\ldots,j_n}(f(\bm{A}_{j_1},\bm{A}_{j_2},\ldots,\bm{A}_{j_n}))}\nonumber \\
&&- g\Bigg(\sum\limits_{j_{1}=1,\ldots,j_{n}=1}^{k_1,\ldots,k_n}w_{j_1}\ldots w_{j_n}\Phi_{j_1,\ldots,j_n}(\bm{A}_{j_1}),\ldots,\sum\limits_{j_{1}=1,\ldots,j_{n}=1}^{k_1,\ldots,k_n}w_{j_1}\ldots w_{j_n}\Phi_{j_1,\ldots,j_n}(\bm{A}_{j_n})\Bigg)
\nonumber \\
&\geq&
\min\limits_{\substack{m_i \leq x_i \leq M_i \\ i=1,2,\ldots,n}}\Big(\sum\limits_{i=1}^n a_i x_i + b - g(x_1,x_2,\ldots,x_n)\Big)\bm{I}_{\mathfrak{K}}.
\end{eqnarray}
\end{theorem}
\textbf{Proof:}
The upper bound of this theorem is proved by setting $\alpha=1$ in Eq.~\eqref{eq UB: thm:main 2.4} from Theorem~\ref{thm:main 2.4} and rearrangement of the term $g\Bigg(\sum\limits_{j_{1}=1,\ldots,j_{n}=1}^{k_1,\ldots,k_n}w_{j_1}\ldots w_{j_n}\Phi_{j_1,\ldots,j_n}(\bm{A}_{j_1}),\ldots,\sum\limits_{j_{1}=1,\ldots,j_{n}=1}^{k_1,\ldots,k_n}w_{j_1}\ldots w_{j_n}\Phi_{j_1,\ldots,j_n}(\bm{A}_{j_n})\Bigg)$ to obtain Eq.~\eqref{eq UB: cor 2.15}. 

Similarly, the lower bound of this theorem is proved by setting $\alpha=1$ in Eq.~\eqref{eq LB: thm:main 2.4} from Theorem~\ref{thm:main 2.4} and rearrangement of the term $g\Bigg(\sum\limits_{j_{1}=1,\ldots,j_{n}=1}^{k_1,\ldots,k_n}w_{j_1}\ldots w_{j_n}\Phi_{j_1,\ldots,j_n}(\bm{A}_{j_1}),\ldots,\sum\limits_{j_{1}=1,\ldots,j_{n}=1}^{k_1,\ldots,k_n}w_{j_1}\ldots w_{j_n}\Phi_{j_1,\ldots,j_n}(\bm{A}_{j_n})\Bigg)$ to obtain Eq.~\eqref{eq LB: cor 2.15}. 
$\hfill \Box$

Next Corollary~\ref{cor: 2.17 ext} is obtained by applying Theorem~\ref{thm: cor 2.15} to special types of function $g(x_1, x_2, \ldots, x_n)$. 
\begin{corollary}\label{cor: 2.17 ext}
Let $\bm{A}_{j_i}$ be self-adjoint operators with $\Lambda(\bm{A}_{j_i}) \in [m_i, M_i]$ for real scalars $m_i <  M_i$. The mappings $\Phi_{j_1,\ldots,j_n}: \mathscr{B}(\mathfrak{H}) \rightarrow \mathscr{B}(\mathfrak{K})$ are normalized positive linear maps, where $j_i=1,2,\ldots,k_i$ for $i=1,2,\ldots,n$. We have $n$ probability vectors $\bm{w}_i =[w_{i,1},w_{i,2},\cdots, w_{i,k_i}]$ with the dimension $k_i$ for $i=1,2,\ldots,n$, i.e., $\sum\limits_{\ell=1}^{k_i}w_{i,\ell} = 1$. Let $f(x_1,x_2,\ldots,x_n)$ be any real valued continuous functions with $n$ variables defined on the range $\bigtimes\limits_{i=1}^n [m_i, M_i] \in \mathbb{R}^n$, where $\times$ is the Cartesian product. Besides, we assume that the function $f$ satisfies the following:
\begin{eqnarray}\label{eq f bounds: cor: 2.17 ext}
\sum\limits_{i=1}^n a_i x_i + b \leq f(x_1,x_2,\ldots,x_n) \leq \sum\limits_{i=1}^n c_i x_i + d,
\end{eqnarray}
for $[x_1,x_2,\ldots,x_n] \in \bigtimes\limits_{i=1}^n [m_i,M_i]$. 

(I) If we also assume that $g(x_1,\ldots,x_n) = \Big(\sum\limits_{i=1}^n \beta_i x_i\Big)^q$, where $\beta_i \geq 0$, $q \in \mathbb{R}$ and $m_i \geq 0$, then, we have the following upper bound for $\sum\limits_{j_{1}=1,\ldots,j_{n}=1}^{k_1,\ldots,k_n}  w_{j_1}\ldots w_{j_n}\Phi_{j_1,\ldots,j_n}(f(\bm{A}_{j_1},\bm{A}_{j_2},\ldots,\bm{A}_{j_n}))$:
\begin{eqnarray}\label{eq x q UB: cor: 2.17 ext}
\lefteqn{\sum\limits_{j_{1}=1,\ldots,j_{n}=1}^{k_1,\ldots,k_n}  w_{j_1}\ldots w_{j_n}\Phi_{j_1,\ldots,j_n}(f(\bm{A}_{j_1},\bm{A}_{j_2},\ldots,\bm{A}_{j_n}))}\nonumber \\
&&- \left(\sum\limits_{i=1}^n \beta_i x_i \sum\limits_{j_{1}=1,\ldots,j_{n}=1}^{k_1,\ldots,k_n}w_{j_1}\ldots w_{j_n}\Phi_{j_1,\ldots,j_n}(\bm{A}_{j_i})\right)^q
\nonumber \\
&\leq&
\max\limits_{\substack{m_i \leq x_i \leq M_i \\ i=1,2,\ldots,n}}\Big(\sum\limits_{i=1}^n c_i x_i + d - \Big(\sum\limits_{i=1}^n \beta_i x_i\Big)^q\Big)\bm{I}_{\mathfrak{K}};
\end{eqnarray}
and, the following lower bound for $\sum\limits_{j_{1}=1,\ldots,j_{n}=1}^{k_1,\ldots,k_n}  w_{j_1}\ldots w_{j_n}\Phi_{j_1,\ldots,j_n}(f(\bm{A}_{j_1},\bm{A}_{j_2},\ldots,\bm{A}_{j_n}))$:
\begin{eqnarray}\label{eq x q LB: cor: 2.17 ext}
\lefteqn{\sum\limits_{j_{1}=1,\ldots,j_{n}=1}^{k_1,\ldots,k_n}  w_{j_1}\ldots w_{j_n}\Phi_{j_1,\ldots,j_n}(f(\bm{A}_{j_1},\bm{A}_{j_2},\ldots,\bm{A}_{j_n}))}\nonumber \\
&&- \left(\sum\limits_{i=1}^n \beta_i x_i \sum\limits_{j_{1}=1,\ldots,j_{n}=1}^{k_1,\ldots,k_n}w_{j_1}\ldots w_{j_n}\Phi_{j_1,\ldots,j_n}(\bm{A}_{j_i})\right)^q
\nonumber \\
&\geq&
\min\limits_{\substack{m_i \leq x_i \leq M_i \\ i=1,2,\ldots,n}}\Big(\sum\limits_{i=1}^n a_i x_i + b - \Big(\sum\limits_{i=1}^n \beta_i x_i\Big)^q\Big)\bm{I}_{\mathfrak{K}}.
\end{eqnarray}

(II) If $g(x_1,\ldots,x_n) = \log\Big(\sum\limits_{i=1}^n \beta_i x_i\Big)$ with $\beta_i \geq 0$ and $m_i >0$, we have the upper bound for $\sum\limits_{j_{1}=1,\ldots,j_{n}=1}^{k_1,\ldots,k_n}  w_{j_1}\ldots w_{j_n}\Phi_{j_1,\ldots,j_n}(f(\bm{A}_{j_1},\bm{A}_{j_2},\ldots,\bm{A}_{j_n}))$:
\begin{eqnarray}\label{eq log x UB: cor: 2.17 ext}
\lefteqn{\sum\limits_{j_{1}=1,\ldots,j_{n}=1}^{k_1,\ldots,k_n}  w_{j_1}\ldots w_{j_n}\Phi_{j_1,\ldots,j_n}(f(\bm{A}_{j_1},\bm{A}_{j_2},\ldots,\bm{A}_{j_n}))}\nonumber \\
&&- \log\left(\sum\limits_{i=1}^n \beta_i x_i \sum\limits_{j_{1}=1,\ldots,j_{n}=1}^{k_1,\ldots,k_n}w_{j_1}\ldots w_{j_n}\Phi_{j_1,\ldots,j_n}(\bm{A}_{j_i})\right)
\nonumber \\
&\leq&
\max\limits_{\substack{m_i \leq x_i \leq M_i \\ i=1,2,\ldots,n}}\Big(\sum\limits_{i=1}^n c_i x_i + d - \log\Big(\sum\limits_{i=1}^n \beta_i x_i\Big)\Big)\bm{I}_{\mathfrak{K}};
\end{eqnarray}
and, the following lower bound for $\sum\limits_{j_{1}=1,\ldots,j_{n}=1}^{k_1,\ldots,k_n}  w_{j_1}\ldots w_{j_n}\Phi_{j_1,\ldots,j_n}(f(\bm{A}_{j_1},\bm{A}_{j_2},\ldots,\bm{A}_{j_n}))$:
\begin{eqnarray}\label{eq log x LB: cor: 2.17 ext}
\lefteqn{\sum\limits_{j_{1}=1,\ldots,j_{n}=1}^{k_1,\ldots,k_n}  w_{j_1}\ldots w_{j_n}\Phi_{j_1,\ldots,j_n}(f(\bm{A}_{j_1},\bm{A}_{j_2},\ldots,\bm{A}_{j_n}))}\nonumber \\
&&- \log\left(\sum\limits_{i=1}^n \beta_i x_i \sum\limits_{j_{1}=1,\ldots,j_{n}=1}^{k_1,\ldots,k_n}w_{j_1}\ldots w_{j_n}\Phi_{j_1,\ldots,j_n}(\bm{A}_{j_i})\right)
\nonumber \\
&\geq&
\min\limits_{\substack{m_i \leq x_i \leq M_i \\ i=1,2,\ldots,n}}\Big(\sum\limits_{i=1}^n a_i x_i + b - \log\Big(\sum\limits_{i=1}^n \beta_i x_i\Big)\Big)\bm{I}_{\mathfrak{K}}.
\end{eqnarray}

(III) If $g(x_1,\ldots,x_n) = \exp\Big(\sum\limits_{i=1}^n \beta_i x_i\Big)$, we have the upper bound for \\
$\sum\limits_{j_{1}=1,\ldots,j_{n}=1}^{k_1,\ldots,k_n}  w_{j_1}\ldots w_{j_n}\Phi_{j_1,\ldots,j_n}(f(\bm{A}_{j_1},\bm{A}_{j_2},\ldots,\bm{A}_{j_n}))$:
\begin{eqnarray}\label{eq exp x UB: cor: 2.17 ext}
\lefteqn{\sum\limits_{j_{1}=1,\ldots,j_{n}=1}^{k_1,\ldots,k_n}  w_{j_1}\ldots w_{j_n}\Phi_{j_1,\ldots,j_n}(f(\bm{A}_{j_1},\bm{A}_{j_2},\ldots,\bm{A}_{j_n}))}\nonumber \\
&&- \exp\left(\sum\limits_{i=1}^n \beta_i x_i \sum\limits_{j_{1}=1,\ldots,j_{n}=1}^{k_1,\ldots,k_n}w_{j_1}\ldots w_{j_n}\Phi_{j_1,\ldots,j_n}(\bm{A}_{j_i})\right)
\nonumber \\
&\leq&
\max\limits_{\substack{m_i \leq x_i \leq M_i \\ i=1,2,\ldots,n}}\Big(\sum\limits_{i=1}^n c_i x_i + d - \exp\Big(\sum\limits_{i=1}^n \beta_i x_i\Big)\Big)\bm{I}_{\mathfrak{K}};
\end{eqnarray}
and, the following lower bound for $\sum\limits_{j_{1}=1,\ldots,j_{n}=1}^{k_1,\ldots,k_n}  w_{j_1}\ldots w_{j_n}\Phi_{j_1,\ldots,j_n}(f(\bm{A}_{j_1},\bm{A}_{j_2},\ldots,\bm{A}_{j_n}))$:
\begin{eqnarray}\label{eq exp x LB: cor: 2.17 ext}
\lefteqn{\sum\limits_{j_{1}=1,\ldots,j_{n}=1}^{k_1,\ldots,k_n}  w_{j_1}\ldots w_{j_n}\Phi_{j_1,\ldots,j_n}(f(\bm{A}_{j_1},\bm{A}_{j_2},\ldots,\bm{A}_{j_n}))}\nonumber \\
&&- \exp\left(\sum\limits_{i=1}^n \beta_i x_i \sum\limits_{j_{1}=1,\ldots,j_{n}=1}^{k_1,\ldots,k_n}w_{j_1}\ldots w_{j_n}\Phi_{j_1,\ldots,j_n}(\bm{A}_{j_i})\right)
\nonumber \\
&\geq&
\min\limits_{\substack{m_i \leq x_i \leq M_i \\ i=1,2,\ldots,n}}\Big(\sum\limits_{i=1}^n a_i x_i + b - \exp\Big(\sum\limits_{i=1}^n \beta_i x_i\Big)\Big)\bm{I}_{\mathfrak{K}}.
\end{eqnarray}
\end{corollary}
\textbf{Proof:}
For Part (I), we will use $g(x_1,\ldots,x_n) = \Big(\sum\limits_{i=1}^n \beta_i x_i\Big)^q$ in Eq.~\eqref{eq UB: cor 2.15} in Theorem~\ref{thm: cor 2.15} to obtain Eq.~\eqref{eq x q UB: cor: 2.17 ext}. We will also use $g(x)=x^q$ in Eq.~\eqref{eq LB: cor 2.15} in Theorem~\ref{thm: cor 2.15} to obtain Eq.~\eqref{eq x q LB: cor: 2.17 ext}. 

For Part (II), we will use  $g(x_1,\ldots,x_n) = \log\Big(\sum\limits_{i=1}^n \beta_i x_i\Big)$ in Eq.~\eqref{eq UB: cor 2.15} in Theorem~\ref{thm: cor 2.15} to obtain Eq.~\eqref{eq log x UB: cor: 2.17 ext}. We will also use $g(x)=\log(x)$ in Eq.~\eqref{eq LB: cor 2.15} in Theorem~\ref{thm: cor 2.15} to obtain Eq.~\eqref{eq log x LB: cor: 2.17 ext}. 

For Part (III), we will use $g(x_1,\ldots,x_n) = \exp\Big(\sum\limits_{i=1}^n \beta_i x_i\Big)$ in Eq.~\eqref{eq UB: cor 2.15} in Theorem~\ref{thm: cor 2.15} to obtain Eq.~\eqref{eq exp x UB: cor: 2.17 ext}. We will also use $g(x)=\exp(x)$ in Eq.~\eqref{eq LB: cor 2.15} in Theorem~\ref{thm: cor 2.15} to obtain Eq.~\eqref{eq exp x LB: cor: 2.17 ext}. 
$\hfill \Box$

\section{Applications}\label{sec: Applications}

In this section, we will apply derived multivariate hypercomplex function inequalities to established Sobolev inequality for hypercomplex function. These inequalities will be used to characterize embedding relationships for multivariate hypercomplex functions. In Section~\ref{sec: Hypercomplex Function Sobolev Embedding By Original Arguments}, Sobolev inequality for hypercomplex function is derived by the original input operators. On the other hand, in Section~\ref{sec: Hypercomplex Function Sobolev Embedding By Arguments Mean}, Sobolev inequality for hypercomplex function is derived based on the arithmetic mean of input operators.

\subsection{Hypercomplex Function Sobolev Embedding By Original Arguments}\label{sec: Hypercomplex Function Sobolev Embedding By Original Arguments}

Before presenting the main result of this subsection, we have to introduce definitions about $|f(\bm{A}_{j_1},\bm{A}_{j_2},\ldots,\bm{A}_{j_n})|$ and $|f'(\bm{A}_{j_1},\bm{A}_{j_2},\ldots,\bm{A}_{j_n})|$. 

We define $|f(\bm{A}_{j_1},\bm{A}_{j_2},\ldots,\bm{A}_{j_n})|$ by
\begin{eqnarray}\label{eq: abs operator def}
|f(\bm{A}_{j_1},\bm{A}_{j_2},\ldots,\bm{A}_{j_n})| \define \sqrt{f^{\ast}(\bm{A}_{j_1},\bm{A}_{j_2},\ldots,\bm{A}_{j_n})f(\bm{A}_{j_1},\bm{A}_{j_2},\ldots,\bm{A}_{j_n})}. 
\end{eqnarray}
Note that $|f(\bm{A}_{j_1},\bm{A}_{j_2},\ldots,\bm{A}_{j_n})|$ is a self-adjoint operator.

We define $|f'(\bm{A}_{j_1},\bm{A}_{j_2},\ldots,\bm{A}_{j_n})|$ by
\begin{eqnarray}\label{eq1: abs der operator def}
|f'(\bm{A}_{j_1},\bm{A}_{j_2},\ldots,\bm{A}_{j_n})| \define \sqrt{\sum\limits_{i=1}^{n}
[f^{(i)}(\bm{A}_{j_1},\bm{A}_{j_2},\ldots,\bm{A}_{j_n})]^2}. 
\end{eqnarray}
where $f^{(i)}(\bm{A}_{j_1},\bm{A}_{j_2},\ldots,\bm{A}_{j_n})$ is expressed by 
\begin{eqnarray}\label{eq2: abs der operator def}
f^{(i)}(\bm{A}_{j_1},\bm{A}_{j_2},\ldots,\bm{A}_{j_n}) = \frac{\partial f(x_1,\ldots,x_n)}{\partial x_i}\Bigg|_{x_1=\bm{A}_{j_1},\ldots,x_n=\bm{A}_{j_n}}.
\end{eqnarray}
Note that $f^{(i)}(\bm{A}_{j_1},\bm{A}_{j_2},\ldots,\bm{A}_{j_n})$ is a self-adjoint operator by our assumption.

We also need Lemma~\ref{lma: C p/q} and Lemma~\ref{lma: C p q} below in order to present our main theorem in this seciton.

\begin{lemma}\label{lma: C p/q}
Given any integer $m > 1$, $p$ be any real number within $1 < p < m$, and $q$ be another real number defined by $q=\frac{mp}{m-p}$. Besides, we are also provided by $n$ self-adjoint operators $\bm{X}_{1},\bm{X}_{2},\ldots,\bm{X}_{n}$ with the spectrum region $\bigtimes_{i=1}^n \Lambda(\bm{X}_{i}) \in \mathbb{R}^n$. If $\Lambda(\bm{X}_{i})$ is a compact region for $i=1,2,\ldots,n$, and the continuous multivaraite function with $h(x_1,x_2,\ldots,x_n)>0$ for $[x_1,x_2,\ldots,x_n] \in \bigtimes_{i=1}^n \Lambda(\bm{X}_{i})$, we can find a positive contant $C_2$ such that    
\begin{eqnarray}\label{eq1: lma: C p/q}
C_{2} h(\bm{X}_{1},\bm{X}_{2},\ldots,\bm{X}_{n}) \leq h^{p/q}(\bm{X}_{1},\bm{X}_{2},\ldots,\bm{X}_{n}).
\end{eqnarray}
\end{lemma}
\textbf{Proof:}
Since $h(x_1,x_2,\ldots,x_n)>0$ for $[x_1,x_2,\ldots,x_n] \in \bigtimes_{i=1}^n \Lambda(\bm{X}_{i})$, we can find a positive $C_2$ such that 
\begin{eqnarray}\label{eq1: lma: C p/q}
C_{2} h(\lambda_1,\lambda_2,\ldots,\lambda_n) \leq h^{p/q}(\lambda_1,\lambda_2,\ldots,\lambda_n),
\end{eqnarray}
where $\lambda_i \in \Lambda(\bm{X}_{i})$, by selection any positive $C_2$ less then the minimum value of $h^{\frac{p-q}{q}}(\lambda_1,\lambda_2,\ldots,\lambda_n)$. Such minimum value exists as a continuous function on a compact set has an absolute minimum. This Lemma follows from spectrum mapping theorem.
$\hfill \Box$

Note that the positive contant $C_2$ depends only on $p$ and $q$ for any multivariate function $h$ bounded from below.

\begin{lemma}\label{lma: C p q}
Given any integer $m > 1$, $p$ be any real number within $1 < p < m$, and $q$ be another real number defined by $q=\frac{mp}{m-p}$. Besides, we are also provided by $n$ self-adjoint operators $\bm{X}_{1},\bm{X}_{2},\ldots,\bm{X}_{n}$ with the spectrum region $\bigtimes_{i=1}^n \Lambda(\bm{X}_{i}) \in \mathbb{R}^n$. If $\Lambda(\bm{X}_{i})$ is a compact region for $i=1,2,\ldots,n$, and the multivaraite function with $|f'(x_1,x_2,\ldots,x_n)|$ bounded below from zero for $[x_1,x_2,\ldots,x_n] \in \bigtimes_{i=1}^n \Lambda(\bm{X}_{i})$, we can find a positive contant $C_3$ such that    
\begin{eqnarray}\label{eq1: lma: C p q}
|f(\lambda_{1},\lambda_{2},\ldots,\lambda_{n})|^p \leq C_3 |f'(\lambda_{1},\lambda_{2},\ldots,\lambda_{n})|^q,
\end{eqnarray}
where $\lambda_i \in \Lambda(\bm{X}_{i})$.
\end{lemma}
\textbf{Proof:}
Because $|f(\lambda_{1},\lambda_{2},\ldots,\lambda_{n})|$ is a continuous function over a compact region, the maximum value of $|f(\lambda_{1},\lambda_{2},\ldots,\lambda_{n})|^p$ is bounded. On the other hand, the multivariate function with $|f'(x_1,x_2,\ldots,x_n)|$ bounded below from zero for $[x_1,x_2,\ldots,x_n] \in \bigtimes_{i=1}^n \Lambda(\bm{X}_{i})$, we also can fine the minimum value of $|f'(\lambda_{1},\lambda_{2},\ldots,\lambda_{n})|^q$, which is a positive number. Therefore, the positive constant $C_3$ can be determined by finding any number greater or equal than the ratio between the maximum value of $|f(\lambda_{1},\lambda_{2},\ldots,\lambda_{n})|^p$ and the minimum value of $|f'(\lambda_{1},\lambda_{2},\ldots,\lambda_{n})|^q$.
$\hfill \Box$

Note that the positive contant $C_3$ depends only on $p$ and $q$ for any multivariate function $f$ satisfying the conditions provided by Lemma~\ref{lma: C p q}.

\begin{theorem}[Sobolev Inequality for Hypercomplex Function]\label{thm: Sobolev Inequality for Hypercomplex Function}
Given any integer $m > 1$, $p$ be any real number within $1 < p < m$, and $q$ be another real number defined by $q=\frac{mp}{m-p}$. Let $\bm{A}_{j_i}$ be self-adjoint operators with $\Lambda(\bm{A}_{j_i}) \in [m_i, M_i]$ for real scalars $m_i <  M_i$. The mappings $\Phi_{j_1,\ldots,j_n}: \mathscr{B}(\mathfrak{H}) \rightarrow \mathscr{B}(\mathfrak{K})$ are normalized positive linear maps, where $j_i=1,2,\ldots,k_i$ for $i=1,2,\ldots,n$. We have $n$ probability vectors $\bm{w}_i =[w_{i,1},w_{i,2},\cdots, w_{i,k_i}]$ with the dimension $k_i$ for $i=1,2,\ldots,n$, i.e., $\sum\limits_{\ell=1}^{k_i}w_{i,\ell} = 1$. If the multivaraite function with $|f'(x_1,x_2,\ldots,x_n)|$ bounded below from zero for $[x_1,x_2,\ldots,x_n] \in \bigtimes_{i=1}^n [m_i, M_i]$. 
We have the following inequality:
\begin{eqnarray}\label{eq1: thm: Sobolev Inequality for Hypercomplex Function}
\lefteqn{\Big(\sum\limits_{j_{1}=1,\ldots,j_{n}=1}^{k_1,\ldots,k_n}  w_{j_1}\ldots w_{j_n}\Phi_{j_1,\ldots,j_n}(|f(\bm{A}_{j_1},\bm{A}_{j_2},\ldots,\bm{A}_{j_n})|^p)\Big)^{1/p}}\nonumber \\
&\leq& C_1\Big(\sum\limits_{j_{1}=1,\ldots,j_{n}=1}^{k_1,\ldots,k_n}  w_{j_1}\ldots w_{j_n}\Phi_{j_1,\ldots,j_n}(|f'(\bm{A}_{j_1},\bm{A}_{j_2},\ldots,\bm{A}_{j_n})|^q)\Big)^{1/q},
\end{eqnarray}
where $C_1$ is a constant indepedent of the function $f(x_1,x_2,\ldots,x_n)$. 

\end{theorem}
\textbf{Proof:}
From Lemma~\ref{lma: C p q} and spectrum mapping theorem, we have
\begin{eqnarray}\label{eq2: thm: Sobolev Inequality for Hypercomplex Function}
\lefteqn{\sum\limits_{j_{1}=1,\ldots,j_{n}=1}^{k_1,\ldots,k_n}  w_{j_1}\ldots w_{j_n}\Phi_{j_1,\ldots,j_n}(|f(\bm{A}_{j_1},\bm{A}_{j_2},\ldots,\bm{A}_{j_n})|^p)}\nonumber \\
&\leq& C_3 \sum\limits_{j_{1}=1,\ldots,j_{n}=1}^{k_1,\ldots,k_n}  w_{j_1}\ldots w_{j_n}\Phi_{j_1,\ldots,j_n}(|f'(\bm{A}_{j_1},\bm{A}_{j_2},\ldots,\bm{A}_{j_n})|^q).
\end{eqnarray}
By rasing the power $1/q$ at both sides of Eq.~\eqref{eq2: thm: Sobolev Inequality for Hypercomplex Function}, which is less than $1$, we have 
\begin{eqnarray}\label{eq3: thm: Sobolev Inequality for Hypercomplex Function}
\lefteqn{\Big(\sum\limits_{j_{1}=1,\ldots,j_{n}=1}^{k_1,\ldots,k_n}  w_{j_1}\ldots w_{j_n}\Phi_{j_1,\ldots,j_n}(|f(\bm{A}_{j_1},\bm{A}_{j_2},\ldots,\bm{A}_{j_n})|^p)\Big)^{1/q}}\nonumber \\
&\leq& C'_3 \Big(\sum\limits_{j_{1}=1,\ldots,j_{n}=1}^{k_1,\ldots,k_n}  w_{j_1}\ldots w_{j_n}\Phi_{j_1,\ldots,j_n}(|f'(\bm{A}_{j_1},\bm{A}_{j_2},\ldots,\bm{A}_{j_n})|^q)\Big)^{1/q},
\end{eqnarray}
where the inequality $\leq$ comes from Lowner-Heinz inequality~\cite{furuta2005mond,fujii2012recentMP}. 

The L.H.S. of Eq.~\eqref{eq3: thm: Sobolev Inequality for Hypercomplex Function} can be expressed by
\begin{eqnarray}\label{eq4: thm: Sobolev Inequality for Hypercomplex Function}
\lefteqn{\Big(\sum\limits_{j_{1}=1,\ldots,j_{n}=1}^{k_1,\ldots,k_n}  w_{j_1}\ldots w_{j_n}\Phi_{j_1,\ldots,j_n}(|f(\bm{A}_{j_1},\bm{A}_{j_2},\ldots,\bm{A}_{j_n})|^p)\Big)^{1/q}}\nonumber \\
&=& \Big(\Big(\sum\limits_{j_{1}=1,\ldots,j_{n}=1}^{k_1,\ldots,k_n}  w_{j_1}\ldots w_{j_n}\Phi_{j_1,\ldots,j_n}(|f(\bm{A}_{j_1},\bm{A}_{j_2},\ldots,\bm{A}_{j_n})|^p)\Big)^{1/p}\Big)^{p/q}
\end{eqnarray}
From Lemma~\ref{lma: C p/q}. we have 
\begin{eqnarray}\label{eq5: thm: Sobolev Inequality for Hypercomplex Function}
\lefteqn{C_2\Big(\sum\limits_{j_{1}=1,\ldots,j_{n}=1}^{k_1,\ldots,k_n}  w_{j_1}\ldots w_{j_n}\Phi_{j_1,\ldots,j_n}(|f(\bm{A}_{j_1},\bm{A}_{j_2},\ldots,\bm{A}_{j_n})|^p)\Big)^{1/p}}\nonumber \\
&\leq& \Big(\Big(\sum\limits_{j_{1}=1,\ldots,j_{n}=1}^{k_1,\ldots,k_n}  w_{j_1}\ldots w_{j_n}\Phi_{j_1,\ldots,j_n}(|f(\bm{A}_{j_1},\bm{A}_{j_2},\ldots,\bm{A}_{j_n})|^p)\Big)^{1/p}\Big)^{p/q}.
\end{eqnarray}
Finally, this Theorem is provoed by combing Eq.~\eqref{eq3: thm: Sobolev Inequality for Hypercomplex Function}, Eq.~\eqref{eq4: thm: Sobolev Inequality for Hypercomplex Function}, Eq.~\eqref{eq5: thm: Sobolev Inequality for Hypercomplex Function}, and setting $C_1 = C'_3/C_2$. 
$\hfill \Box$

We define $\mathscr{W}_{\Phi}^{1,q}\left(\bigtimes_{i=1}^n [m_i, M_i]\right)$ as the hypercomplex function space such that 
\begin{eqnarray}
\sum\limits_{j_{1}=1,\ldots,j_{n}=1}^{k_1,\ldots,k_n}  w_{j_1}\ldots w_{j_n}\Phi(|f'(\bm{A}_{j_1},\bm{A}_{j_2},\ldots,\bm{A}_{j_n})|^q)\Big)^{1/q} 
\end{eqnarray}
is a bounded and positive self-adjoint operator for spectrum in the range of $\bigtimes_{i=1}^n [m_i, M_i]$. In addition, we define $\mathscr{L}^{p}_{\Phi}\left(\bigtimes_{i=1}^n [m_i, M_i]\right)$
as the hypercomplex function space such that 
\begin{eqnarray}
\sum\limits_{j_{1}=1,\ldots,j_{n}=1}^{k_1,\ldots,k_n}  w_{j_1}\ldots w_{j_n}\Phi(|f(\bm{A}_{j_1},\bm{A}_{j_2},\ldots,\bm{A}_{j_n})|^p)\Big)^{1/p}
\end{eqnarray}
is a bounded and positive self-adjoint operator for spectrum in the range of $\bigtimes_{i=1}^n [m_i, M_i]$. From Theorem~\ref{thm: Sobolev Inequality for Hypercomplex Function}, we have hypercomplex function Sobolev embedding as 
\begin{eqnarray}
\mathscr{W}_{\Phi}^{1,q}\left(\bigtimes_{i=1}^n [m_i, M_i]\right)\subset \mathscr{L}^{p}_{\Phi}\left(\bigtimes_{i=1}^n [m_i, M_i]\right).
\end{eqnarray}

\subsection{Hypercomplex Function Sobolev Embedding By Arguments Mean}\label{sec: Hypercomplex Function Sobolev Embedding By Arguments Mean}

\begin{theorem}[Sobolev Inequality for Hypercomplex Function with Operator Mean]\label{thm: Sobolev Inequality for Hypercomplex Function mean}
Given any integer \\
$m > 1$, $p$ be any real number within $1 < p < m$, and $q$ be another real number defined by $q=\frac{mp}{m-p}$. Let $\bm{A}_{j_i}$ be self-adjoint operators with $\Lambda(\bm{A}_{j_i}) \in [m_i, M_i]$ for real scalars $m_i <  M_i$. The mappings $\Phi_{j_1,\ldots,j_n}: \mathscr{B}(\mathfrak{H}) \rightarrow \mathscr{B}(\mathfrak{K})$ are normalized positive linear maps, where $j_i=1,2,\ldots,k_i$ for $i=1,2,\ldots,n$. We have $n$ probability vectors $\bm{w}_i =[w_{i,1},w_{i,2},\cdots, w_{i,k_i}]$ with the dimension $k_i$ for $i=1,2,\ldots,n$, i.e., $\sum\limits_{\ell=1}^{k_i}w_{i,\ell} = 1$. 

Besides, we assume that the function $f$ satisfies the following:
\begin{eqnarray}\label{eq f bounds: thm: Sobolev Inequality for Hypercomplex Function mean}
|f(x_1,x_2,\ldots,x_n)|^p \leq \sum\limits_{i=1}^n c_i x_i + d,
\end{eqnarray}
where $[x_1,x_2,\ldots,x_n] \in \bigtimes_{i=1}^n [m_i, M_i]$. In addition, we also assume that the function $|f'|^q$ satisfies the following:
\begin{eqnarray}\label{eq g bounds: thm: Sobolev Inequality for Hypercomplex Function mean}
|f'(x_1,x_2,\ldots,x_n)|^q > 0,
\end{eqnarray}
where $[x_1,x_2,\ldots,x_n] \in \bigtimes_{i=1}^n [m_i, M_i]$. We have the following inequality:
\begin{eqnarray}\label{eq1: thm: Sobolev Inequality for Hypercomplex Function mean}
\lefteqn{\Bigg(\sum\limits_{j_{1}=1,\ldots,j_{n}=1}^{k_1,\ldots,k_n}  w_{j_1}\ldots w_{j_n}\Phi_{j_1,\ldots,j_n}(|f(\bm{A}_{j_1},\bm{A}_{j_2},\ldots,\bm{A}_{j_n})|^p)\Bigg)^{1/p}\leq}\nonumber \\
& &C_4 \Bigg(\Bigg|f'\Big(\sum\limits_{j_{1}=1,\ldots,j_{n}=1}^{k_1,\ldots,k_n}w_{j_1}\ldots w_{j_n}\Phi_{j_1,\ldots,j_n}(\bm{A}_{j_1}),\ldots,\sum\limits_{j_{1}=1,\ldots,j_{n}=1}^{k_1,\ldots,k_n}w_{j_1}\ldots w_{j_n}\Phi_{j_1,\ldots,j_n}(\bm{A}_{j_n})\Bigg)\Bigg|^q\Bigg)^{1/q}, \nonumber \\
\end{eqnarray}
where $C_4$ is a positive scalar related to $f$, $f'$, $c_i$, $d$, and the domain $\bigtimes_{i=1}^n [m_i, M_i]$. 
\end{theorem}
\textbf{Proof:}
By setting $f=|f(x_1,x_2,\ldots,x_n)|^p$ and $g=|f'(x_1,x_2,\ldots,x_n)|^q$ in Part (I) of Theorem~\ref{thm: 2.9}, we have 
\begin{eqnarray}\label{eq2: thm: Sobolev Inequality for Hypercomplex Function mean}
\lefteqn{\sum\limits_{j_{1}=1,\ldots,j_{n}=1}^{k_1,\ldots,k_n}  w_{j_1}\ldots w_{j_n}\Phi_{j_1,\ldots,j_n}(|f(\bm{A}_{j_1},\bm{A}_{j_2},\ldots,\bm{A}_{j_n})|^p)}\nonumber \\
&\leq&\left[\max\limits_{\substack{m_i \leq x_i \leq M_i \\ i=1,2,\ldots,n}}\frac{\sum\limits_{i=1}^n c_i x_i +d}{|f'(x_1,x_2,\ldots,x_n)|^q}\right]\nonumber \\
&&\times \Bigg|f'\Bigg(\sum\limits_{j_{1}=1,\ldots,j_{n}=1}^{k_1,\ldots,k_n}w_{j_1}\ldots w_{j_n}\Phi_{j_1,\ldots,j_n}(\bm{A}_{j_1}),\ldots,\sum\limits_{j_{1}=1,\ldots,j_{n}=1}^{k_1,\ldots,k_n}w_{j_1}\ldots w_{j_n}\Phi_{j_1,\ldots,j_n}(\bm{A}_{j_n})\Bigg)\Bigg|^q.
\end{eqnarray}
By rasing the power $1/p$ at both sides of Eq.~\eqref{eq2: thm: Sobolev Inequality for Hypercomplex Function mean}, which is less than $1$, we have 
\begin{eqnarray}\label{eq3: thm: Sobolev Inequality for Hypercomplex Function mean}
\lefteqn{\Bigg(\sum\limits_{j_{1}=1,\ldots,j_{n}=1}^{k_1,\ldots,k_n}  w_{j_1}\ldots w_{j_n}\Phi_{j_1,\ldots,j_n}(|f(\bm{A}_{j_1},\bm{A}_{j_2},\ldots,\bm{A}_{j_n})|^p)\Bigg)^{1/p}}\nonumber \\
&\leq_1&\Bigg(\left[\max\limits_{\substack{m_i \leq x_i \leq M_i \\ i=1,2,\ldots,n}}\frac{\sum\limits_{i=1}^n c_i x_i +d}{|f'(x_1,x_2,\ldots,x_n)|^q}\right]\nonumber \\
&&\times \Bigg|f'\Bigg(\sum\limits_{j_{1}=1,\ldots,j_{n}=1}^{k_1,\ldots,k_n}w_{j_1}\ldots w_{j_n}\Phi_{j_1,\ldots,j_n}(\bm{A}_{j_1}),\ldots,\sum\limits_{j_{1}=1,\ldots,j_{n}=1}^{k_1,\ldots,k_n}w_{j_1}\ldots w_{j_n}\Phi_{j_1,\ldots,j_n}(\bm{A}_{j_n})\Bigg)\Bigg|^q\Bigg)^{1/p}\nonumber \\
&\leq_2&C'_4\Bigg(\left[\max\limits_{\substack{m_i \leq x_i \leq M_i \\ i=1,2,\ldots,n}}\frac{\sum\limits_{i=1}^n c_i x_i +d}{|f'(x_1,x_2,\ldots,x_n)|^q}\right]\nonumber \\
&&\times \Bigg|f'\Bigg(\sum\limits_{j_{1}=1,\ldots,j_{n}=1}^{k_1,\ldots,k_n}w_{j_1}\ldots w_{j_n}\Phi_{j_1,\ldots,j_n}(\bm{A}_{j_1}),\ldots,\sum\limits_{j_{1}=1,\ldots,j_{n}=1}^{k_1,\ldots,k_n}w_{j_1}\ldots w_{j_n}\Phi_{j_1,\ldots,j_n}(\bm{A}_{j_n})\Bigg)\Bigg|^q\Bigg)^{1/q},
\end{eqnarray}
where the inequality $\leq_1$ comes from Lowner-Heinz inequality, and  the inequality $\leq_2$ comes from a positive scalar such that $C'_4 \lambda^{p/q} > \lambda$ for any eigenvalues $\lambda$ of the following term:
\begin{eqnarray}
\Bigg|f'\Bigg(\sum\limits_{j_{1}=1,\ldots,j_{n}=1}^{k_1,\ldots,k_n}w_{j_1}\ldots w_{j_n}\Phi_{j_1,\ldots,j_n}(\bm{A}_{j_1}),\ldots,\sum\limits_{j_{1}=1,\ldots,j_{n}=1}^{k_1,\ldots,k_n}w_{j_1}\ldots w_{j_n}\Phi_{j_1,\ldots,j_n}(\bm{A}_{j_n})\Bigg)\Bigg|^q.
\end{eqnarray}
Finally, this Theorem is proved by setting $C_4$ as 
\begin{eqnarray}
C_4 = C'_4\left[\max\limits_{\substack{m_i \leq x_i \leq M_i \\ i=1,2,\ldots,n}}\frac{\sum\limits_{i=1}^n c_i x_i +d}{|f'(x_1,x_2,\ldots,x_n)|^q}\right]^{1/q}.
\end{eqnarray}
$\hfill \Box$

We define $\overline{\mathscr{W}}_{\Phi}^{1,q}\left(\bigtimes_{i=1}^n [m_i, M_i]\right)$ as the hypercomplex function space  in terms of operator mean such that 
\begin{eqnarray}
\Bigg(\Bigg|f'\Big(\sum\limits_{j_{1}=1,\ldots,j_{n}=1}^{k_1,\ldots,k_n}w_{j_1}\ldots w_{j_n}\Phi(\bm{A}_{j_1}),\ldots,\sum\limits_{j_{1}=1,\ldots,j_{n}=1}^{k_1,\ldots,k_n}w_{j_1}\ldots w_{j_n}\Phi(\bm{A}_{j_n})\Bigg)\Bigg|^q\Bigg)^{1/q}
\end{eqnarray}
is a bounded and positive self-adjoint operator for spectrum in the range of $\bigtimes_{i=1}^n [m_i, M_i]$. From Theorem~\ref{thm: Sobolev Inequality for Hypercomplex Function mean}, we have hypercomplex function Sobolev embedding in terms of operator mean as 
\begin{eqnarray}
\overline{\mathscr{W}}_{\Phi}^{1,q}\left(\bigtimes_{i=1}^n [m_i, M_i]\right)\subset \mathscr{L}^{p}_{\Phi}\left(\bigtimes_{i=1}^n [m_i, M_i]\right).
\end{eqnarray}

\bibliographystyle{IEEETran}
\bibliography{MultivariateMP_Linear_Bib}

\end{document}